\documentclass[11pt]{amsart}

\usepackage{times}
\usepackage{amssymb}
\usepackage{hyperref}

\newtheorem{theo}{Theorem}[section]
\newtheorem{prop}[theo]{Proposition}
\newtheorem{lemm}[theo]{Lemma}
\newtheorem{defi}[theo]{Definition}

\newtheorem{coro}[theo]{Corollary}

\newcommand{\demo}{\begin{proof}}
\newcommand{\findemo}{\end{proof}}

\newcommand{\ba}[1]{\begin{array}{#1}}
\newcommand{\ea}{\end{array}}
\newcommand{\be}{\begin{eqnarray}}
\newcommand{\ee}{\end{eqnarray}}
\newcommand{\benn}{\begin{eqnarray*}}
\newcommand{\eenn}{\end{eqnarray*}}
\newcommand{\norm}[1]{\left\|{#1}\right\|}

\newcommand{\grp}[2]{{#1}\rightrightarrows{#2}}
\newcommand{\op}[3]{C^\infty_c({#1})\otimes_{#2}C^\infty_c({#3})}

\begin{document}

\def\C{\mathbb{C}}
\def\R{\mathbb{R}}
\def\N{\mathbb{N}}
\def\Z{\mathbb{Z}}
\def\T{\mbox{T}}
\def\lb{\langle}
\def\rb{\rangle}
\def\del{\partial}
\def\d{\mbox{d}}
\def\lra{\longrightarrow}
\def\lmt{\longmapsto}
\def\H{\mbox{\bf H}}
\def\V{\mbox{\bf V}}
\def\D{\mbox{D}}
\def\exp{Exp}
\def\pr{\mbox{pr}}
\def\ker{\mbox{Ker}}
\def\im{\mbox{Im}}
\def\id{\mbox{Id}}
\def\S{\mathfrak{S}}
\def\G{\mathfrak{G}}
\def\F{\mathfrak{F}}
\def\D{\mathfrak{D}}
\def\Th{\mathcal{T}}
\def\Dv{\D_v}
\def\half{\frac{1}{2}}
\def\inv{^{-1}}
\def\supp{\mbox{supp}}
\def\sup{\mbox{sup}}
\def\difh{\frac{\d}{\d\hbar}}
\def\difhz{\left.\difh\right|_{\hbar=0}}
\def\cl{\mathnormal{Cl}}

\title{Quantisation of Lie-Poisson manifolds}
\author{S\'ebastien Racani\`ere}
\email{s.racaniere@ic.ac.uk}
\date{\today}

\begin{abstract}
In quantum physics, the operators associated with the position and the momentum of a
particle are unbounded operators and $C^*$-algebraic quantisation does therefore not
deal with such operators. In the present article, I propose a quantisation of the Lie-Poisson
structure of the dual of a Lie algebroid which deals with a big enough class of functions
to include the above mentioned example.

As an application, I show with an example
how the quantisation of the dual of the Lie algebroid associated to a Poisson manifold
can lead to a quantisation of the Poisson manifold itself. The example I consider is
the torus with constant Poisson structure, in which case I recover its usual
$C^*$-algebraic quantisation.
\end{abstract}

\maketitle
\tableofcontents

\section*{Acknowledgement}
I would like to thank B.~Ramazan for sending me a copy of his PhD thesis.

The author was supported by a Marie Curie Fellowship EC Contract Ref:
HPMCF-CT-2002-01850.

\section{Introduction}
In his PhD thesis, B.Ramazan~\cite{Ra} (see also N.~P.~Landsman and B.~Ramazan~\cite{LR})
proved a conjecture of Landsman which roughly
speaking states that the quantised, that is deformed, algebra of functions on the dual of
a Lie algebroid in the direction of its natural Lie-Poisson bracket is the $C^*$-algebra
of the Lie groupoid integrating the Lie algebroid\footnote{That is
if the Lie algebroid is integrable.} The type of quantisations that Ramazan considers are
deformation quantisations in the sense of M.~A.~Rieffel~\cite{Rie}. Not all functions are
quantised in this way, in
fact only functions whose Fourier transform is compactly supported
(with respect to a given family of measures) are quantised.

If $M$ is a Riemannian manifold, then its tangent bundle is a Lie algebroid which integrates
to the pair groupoid $M\times M$. The induced Lie-Poisson structure on $T^*M$ is the usual
symplectic structure on a cotangent bundle and Landsman-Ramazan's quantisation can
be carried over. Nevertheless, this example shows an important limitation of this procedure:
functions that are polynomials in the fibres of $T^*M\lra M$ are not quantised, whereas functions
giving the position or the momentum of a particle are of this type. Moreover, it is well-known
by physicists that the quantisation of such functions are unbounded operators, whereas
Landsman-Ramazan's quantisation only gives elements of $C^*$-algebras, that is bounded operators
on a Hilbert space.

In the present article, I wish to propose a quantisation of the dual of an integrable Lie
algebroids $A\lra M$ which can be used on a wide class of functions. This class contains in
particular functions which are polynomial in the fibres of $A^*\lra M$. This will be done in
Sections~\ref{sect:Method-of-quantisation} and \ref{sect:Computation-local-coordinates}, where
Theorem~\ref{theo:general-quantisation} is the main result. In
Section~\ref{sect:quantisation-of-Rn}, I show that Theorem~\ref{theo:general-quantisation}
can be used to recover the physicists's position and momentum operators of a particle moving
in $\R^n$.

If $M$ is a Poisson manifold, then its cotangent bundle is naturally a Lie algebroid whose
dual can be quantised using Theorem~\ref{theo:general-quantisation}. One might then hope that
this quantisation will help finding a quantisation of the original Poisson manifold $M$. This
slightly naive idea is shown to work on an example, the torus with constant Poisson structure,
in Section~\ref{sect:general-results-about-poisson} and Section~\ref{sect:the-torus}. There,
I recover the usual
$C^*$-algebraic quantisation of a constant Poisson structure on a torus (see
X.~Tang and A.~Weinstein~\cite{TW}, and A.~Weinstein~\cite{Wei}). Part of the strategy of
Section~\ref{sect:the-torus} consists in finding a Poisson map between $T M$ and $M$. Such maps
are solutions to a partial differential equation derived in
Section~\ref{sect:general-results-about-poisson}. In Appendix~A, I show how to find a solution to
this equation in the case of the sphere in $\R^3$.


\section{Method of quantisation}\label{sect:Method-of-quantisation}
Let $\grp{G}{M}$ be a groupoid. Denote by $s$ and 
$t$ the respective source and target maps of the groupoid $\grp{G}{M}$. Also,
denote by $G^{(2)}=\{(x,y)\in G^2\mid\,s(x)=t(y)\}$ the set of pairs of composable arrows in $G$.
Let $\tau:A\lra M$ be Lie algebroid of $\grp{G}{M}$. I will use the same
letter $\tau$ to denote the projection $A^*\lra M$ of the dual of $A$. Choose a Riemannian
metric on $A\lra M$. By duality, this also gives a Riemannian metric on $A^*$. I will denote
by $X$, $Y$ or $Z$ elements in $A^*$ and by $\xi$ or $\zeta$ elements in $A$.

\begin{defi}
Let $E$ be a $s$-family of operators on $G$, that is a map
$q\lmt E_q$ from $M$ to the linear forms on $C^\infty_c(s\inv(q))$. 
I will denote $\op{G}{s}{M}$ the vector space of such operators which in addition satisfy:
for all
smooth family of functions $H$ on $G$ with compact support, that is for all compactly
supported smooth function $H$ on $N\times G$ for some manifold $N$, the function
$$\ba{ccc}
N\times M & \lra & \C\\
(u,q) & \lmt & E_q(x\lmt H(u,x))
\ea$$
is smooth and compactly supported.

Also I will denote $\mathfrak{Op}(G)$ the vector space of $s$-family of operators $E$
which in addition satisfy: for all
smooth family of functions $H$ on $G$ with compact support, that is for all compactly
supported smooth function $H$ on $N\times G$ for some manifold $N$, the function
$$\ba{ccc}
N\times G & \lra & \C\\
(u,z) & \lmt & E_{t(z)}(x\lmt H(u,x z))
\ea$$
is smooth and compactly supported.
\end{defi}

Notice that since $M$ is a closed sub-manifold of $G$, the space $\mathfrak{Op}(G)$ is
included in $\op{G}{s}{M}$.
On the contrary,
\begin{prop}
Let $D$ be in $\op{G}{s}{M}$.
For any compactly supported smooth function $H$ on $N\times M$, the map
$$\ba{ccc}
N\times G & \lra & \C\\
(u,z) & \lmt & D_{t(z)}(H(u,\cdot z))
\ea$$
is smooth. Nevertheless, it might fail to be compactly supported.
\end{prop}
\begin{proof}
Let $H$ be a compactly supported smooth function on $N\times G$, where $N$ is a manifold.
I wish to prove that the map
$$\ba{ccc}
N\times G & \lra & \C\\
(u,z) & \lmt & D_{t(z)}(H(u,\cdot z)))
\ea$$
is smooth.
Let $(u_0,z_0)$ be a point in $N\times G$.
Let $\varphi$ be a compactly supported smooth function on $G$ such that $\varphi\equiv 1$ on
a neighbourhood of $z_0$.
Consider
$$\ba{cccc}
\widetilde{H}: & N\times G^{(2)} & \lra & \C\\
 & (u,x,z) & \lmt & H(u,x\cdot z)\varphi(z).
\ea$$
If
$$K^\prime=\{(u,x,z)\in N\times G^{(2)}\mid\, (u,xz)\in\supp H,\,z\in\supp\varphi\},$$
then the support of $\widetilde{H}$ is a closed subset of $K^\prime$ and since one
easily checks that $K^\prime$ is compact, it follows that $\widetilde{H}$ has compact
support. Using the fact that $G^{(2)}$ is a closed sub-manifold of $G\times G$,
I extend $\widetilde{H}$ to a function, still denoted $\widetilde{H}$, in
$C^\infty_c(N\times G\times G)$. Let $\widetilde{N}=N\times G$. By interpreting the
extended version of $\widetilde{H}$ as a function
$$\ba{ccc}
\widetilde{N}\times G & \lra & \C\\
((u,z),x) & \lmt & \widetilde{H}(u,x,z),
\ea$$
I can apply $D$ and obtain a function
$$\ba{ccc}
\widetilde{N}\times M & \lra & \C\\
((u,z),q) & \lmt & D_q(\widetilde{H}(u,\cdot,z))
\ea$$
in $C^\infty_c(\widetilde{N}\times M)$. The closed sub-manifold $\{x,z,t(z)\}$ of
$\widetilde{N}\times G$ is diffeomorphic to $N\times G$. Therefore
$(u,z)\lmt D_{t(z)}(H(u,\cdot z)\varphi(z))$ is smooth. This map is equal to
$(u,z)\lmt D_{t(z)}(H(u,\cdot z))$ in a neighbourhood of $(u_0,z_0)$. This can be done
for any choice of $(u_0,z_0)$, it follows that $(u,z)\lmt D_{t(z)}(H(u,\cdot z))$ is
smooth.

Nevertheless, the map $(u,z)\lmt D_{t(z)}(H(u,\cdot z))$ needs not be compactly
supported. Indeed, choose $M$ to be a point, that is $G$ is a genuine group. Let $f$ be
a smooth function on $G$, $\mu=\d g$ be a right invariant measure on $G$
and $D=D_f$ be defined as in Proposition~\ref{prop:Definition-of-Df}. Then for $h$ a function on $G$, the map
$z\lmt\int_G f(g)h(g z)\d g$ is certainly not compactly supported in general. For example,
if $f$ is the constant function equal to $1$ then the above map is the constant function
equal to $\int_G h\d g$.
\end{proof}

Recall that $A\lra M$ is the Lie algebroid of $\grp{G}{M}$. Let $T^s G\lra G$ be the
vector bundle whose fibre above $x$ in $G$ is $\ker s_{*,x}$. Denote by
$\mid\Omega\mid^1(T^s G)$ the vector bundle of $1$-densities on the fibres
of $T^s G\lra G$.
Assume that we have a right invariant everywhere positive section $\mu$ of
$\mid\Omega\mid^1(T^s G)$; it defines a  right invariant
smooth Haar system on $\grp{G}{M}$. This section is entirely determined by its value
along $M$ in $G$, which is a section, denoted by $\d\mu$, of $\mid\Omega\mid^1(A)$.
Equivalently, $\d\mu$ is a smooth family of Lebesgue measures on the fibres of $A\lra M$.
By taking the dual, we obtain a family of Lebesgue
measures on $A^*\lra M$, the dual of $A\lra M$.

Integration provides a way of embedding $C^\infty(G)$ in $\op{G}{s}{M}$.

\begin{prop}\label{prop:Definition-of-Df}
Let $f$ be in $C^\infty(G)$.
For each $q$ in $M$, consider the following linear form
on $C^\infty_c(s\inv(q))$
$$D_{f,q}:h\lmt \int_{s\inv(q)}f h\mu.$$
Then $D_f$ is in $\op{G}{s}{M}$.

If moreover $f$ has compact support, then $D_f$ is in $\mathfrak{Op}(G)$.
\end{prop}
\begin{proof}
Let $N$ be a manifold and $H$ a compactly supported smooth function on $N\times G$. It is
clear that the map $(u,q)\lmt D_{f,q}(H(u))$ has support in $(\id_N\times s)(\supp_{H})$
which is compact.

Moreover, because $f H$ has compact support, it is a finite sum of functions with support
contained in open local charts. Writing things in
these local coordinates, it becomes obvious that the map $(u,q)\lmt D_{f,q}(H(u))$ is
smooth.

In addition, when $f$ has compact support, the map
\benn
N\times G^{(2)} & \lra & \C\\
(u,x,z) & \lmt & f(x)H(u,xz)
\eenn
has compact support. Hence
\benn
N\times G & \lra & \C\\
(u,z) & \lmt & \int_{s\inv(t(z))}f(x)H(u,xz)
\eenn
has compact support and $D_f$ is in $\mathfrak{Op}(G)$.
\end{proof}

The Lie algebroid $\tau:A\lra M$ is in particular a vector bundle and one can construct a
Lie groupoid $\grp{A}{M}$ with both the source
and the target map equal to the projection $\tau:A\lra M$. In particular, each smooth
function on $A$ gives an element of $\op{A}{\tau}{M}$.
Let $f$ be a smooth function on $A^*$. Unless the restriction of $f$ to each fibre of $A^*$
is $L^1$, the Fourier transform of $f$ is not defined. Nevertheless, the Fourier transform of
$D_f$ is defined for a much larger class of functions.
\begin{defi}
Let $f$ be a smooth function on $A^*$. Say that $f$ has {\bf polynomial controlled growth} if
\begin{itemize}
\item for every $q$ in $M$,
\item every smooth multi-vector field $\upsilon$ on $M$,
\item every non-negative integer $k$ and every section $\delta$ of $S^k A^*$, and
\item every trivialisation $\phi:A|_{B^\prime}\lra A_q\times B^\prime$ in a neighbourhood
$B^\prime$ of $q$,
\end{itemize}
there exists a smaller neighbourhood $B\subset B^\prime$ of $q$, a non negative constant $C$
and an integer $m$ such that
\be
(\Upsilon\Delta\cdot f)(Y)\leq C(1+\|Y\|^2)^m,\,\mbox{for all }
Y\mbox{ in }A^*|_B,\label{equa:controlled-polynomial-growth}
\ee
where $\Upsilon$ is the multi-vector field defined on $A^*$ using $\upsilon$ and the
trivialisation $\phi$, and $\Delta$ is the multi-vector field on $A^*$ defined using
$\delta$ and the vector space structure on the fibres of $A^*$.

Denote by $C^\infty_{pg}(A^*)$ the set of smooth functions on $A^*$ with
polynomial controlled growth.
\end{defi}
Notice that the above Definition remains unchanged if one replaces
$\Upsilon\Delta$ by $\Delta\Upsilon$ in (\ref{equa:controlled-polynomial-growth}).
Also, to check if $f$ has polynomial controlled growth, it is enough to check
(\ref{equa:controlled-polynomial-growth}) for only one particular choice of trivialisation
$\phi$.

The interesting thing about functions with polynomial controlled growth is that one 
can define the Fourier transform of the operator $D_f$ associated to them. This will be a
consequence of the following easy Lemma.
\begin{lemm}\label{lemm:Fourier-transform-of-compactly-supp-is-Schwartz}
Let $(t,\xi)$ be coordinates on $\R\times\R^n$ and $K$ a compactly supported smooth function
on $\R\times\R^n$. If $P$ is any polynomial function on $\R^n$ then the map
\benn
\R\times\R^n & \lra & \C\\
(t,X) & \lmt & P(X)\int_{\R^n}\d\xi\,e^{-i\langle X,\xi\rangle}K(t,\xi)
\eenn
is bounded.
\end{lemm}
\begin{proof}
This is just a simple application of the fact that the Fourier transform takes
multiplication by a variable to differentiation with respect to that variable.
\end{proof}

If $h$ is a $L^1$ function on a vector space $E$ with measure $\d(X)$, set its Fourier
transform to be the following function on the dual $E^*$ of $E$
\benn
\mathfrak{F}(h)(\xi)=\int_E\d(X)\;e^{-i\lb\xi,X\rb}h(X).
\eenn

\begin{coro}
If $f$ has polynomial controlled growth, set the Fourier transform of $D_f$ to be
\benn
\mathfrak{F}(D_f)_q(h)=\int_{A^*_q}\d\mu(X) f(X)\mathfrak{F}(h)(X),\,\forall q\in M,\,
h\in C^\infty_c(\tau\inv(q)).
\eenn
This is a well-defined element of $\op{A}{\tau}{M}$.
\end{coro}
\begin{proof}
If $H$ is a compactly supported smooth function on $A$, then using
Lemma~\ref{lemm:Fourier-transform-of-compactly-supp-is-Schwartz}
it is easy to prove that $D_f(H)$ is a well-defined smooth function on $Q$.

In addition, the support of $D_f(H)$ is included in the image of the support of $H$ under
the projection $\tau:A\lra Q$ and is therefore compact.
\end{proof}

In Definition~\ref{defi:When-Is-f-Acceptable-ForQuantisation},
I define the set of functions acceptable for quantisation as
a subset of the set of functions with polynomial controlled growth. One can then apply
Proposition~\ref{prop:a-big-set-of-functions-acceptable-for-quantisation}
to see that the set of functions with polynomial controlled growth is
big enough for our purpose. Moreover it is a Poisson algebra as the next Lemma
shows.

\begin{lemm}
The set of functions with polynomial controlled growth forms a Poisson sub-algebra
of $C^\infty(A^*)$.
\end{lemm}
\begin{proof}
This is a simple consequence of Lemma~\ref{lemm:poisson-bracket-of-f-and-g}.
\end{proof}

I now wish to put a structure of algebra on $\mathfrak{Op}(G)$.
\begin{prop}
Let $D$ and $E$ be two elements of $\mathfrak{Op}(G)$.
For $q$ in $M$ and $h$ in $C^\infty_c(s\inv(q))$, set
$$(D\star E)_q(h)=E_q(z\lmt D_{t(z)}(R^*_z h)).$$
The operator $D\star E$ lies in $\mathfrak{Op}(G)$.
\end{prop}
\begin{proof}
Let $N$ be a manifold and $H$ a compactly supported smooth function on $N\times G$.
The function $F$ on $N\times G$ defined by $F(u,z)=D_{t(z)}(H(u,\cdot z))$ is smooth
and compactly supported because $D$ is in $\mathfrak{Op}(G)$; therefore the function
$(u,z)\lmt E_{t(z)}(F(u,\cdot z))$ is smooth and compactly supported.
\end{proof}

Let $\widetilde{G}$ be the tangent groupoid of $G$ (see M.~Hilsum and
G.~Skandalis~\cite{HS}, A.Weinstein~\cite{Wei2} and~\cite{Wei3}, and A.~Connes~\cite{Co})
with respective source and
target maps $\widetilde{s}$ and $\widetilde{t}$. Let
$\widetilde{\tau}:\widetilde{A}=\R\times A\lra \R\times M$ be the
Lie algebroid of the tangent groupoid (this is the tangent algebroid, see
V.~Nistor, A.~Weinstein and P.~Xu~\cite{NWX}).
I will present a method to construct a map from $\op{A}{\tau}{M}$ to
$\op{\widetilde{G}}{\widetilde{s}}{\R\times M}$.
Let $\alpha$ be a diffeomorphism from an open neighbourhood $W$ of $M$ in $A$ to
an open neighbourhood $V$ of $M$ in $G$ such that
\begin{itemize}
\item $\alpha(q)=q$ for $q$ in $M$,
\item $s\circ\alpha=\tau$, in particular $\alpha$ sends $A_q$ to $s\inv(q)$,
\item the differential at zero of the restriction of $\alpha$ to $A_q$ is
the identity map from $A_q$ to $A_q$.
\end{itemize}
For example such an $\alpha$ can be obtained from the choice of an exponential map.
Let
$$\widetilde{W}=\{(\hbar,X)\in\R\times A\mid\, \hbar X\in W\}$$
be an open subset in $\widetilde{A}$. On it, the map
$$\widetilde{\alpha}(\hbar,X)=\left\{\ba{l}
(\hbar,\alpha(\hbar X))\mbox{ for }\hbar\neq 0\\
(0,X)\mbox{ for }\hbar=0
\ea\right.$$
is a diffeomorphism onto an open neighbourhood $\widetilde{V}$ of $\R\times M$ in $\widetilde{G}$.
Choose a smooth function $\psi$ on $A$ with support in $W$ such
that $\psi\mid_{\half W}\equiv 1$. Define $\widetilde{\psi}$ in $C^\infty(\R\times A)$
by $\widetilde{\psi}(\hbar,X)=\psi(\hbar X)$.

\begin{prop}\label{prop:definition-of-tilde}
Let $D$ be in $\op{A}{\tau}{M}$. For $(\hbar,q)$ in $\R\times M$ and
$h$ in $C^\infty_c(\widetilde{s}\inv(\hbar,q))$, set
$$\widetilde{D}_{(\hbar,q)}(h)=D_q(X\lmt \widetilde{\psi}(\hbar,X)h\circ\widetilde{\alpha}(\hbar,X)).$$
The operator $\widetilde{D}$ lies in $\op{\widetilde{G}}{\widetilde{s}}{\R\times M}$.
\end{prop}
\begin{proof}
Let $N$ be a smooth manifold and $H$ a compactly supported smooth function on
$N\times \widetilde{G}$. The function $\widetilde{\psi}\circ\widetilde{\alpha}\inv$
defined on $\widetilde{V}$ can be extended, by zero, to a smooth function on the whole
of $\widetilde{G}$. The product of this function with $H$ is of course with compact
support in $\widetilde{V}$; hence its pull back by $\widetilde{\alpha}$ is compactly
supported. It follows that the function
$$\ba{ccc}
(\N\times\R)\times A & \lra & \mathbb{C}\\
((u,\hbar),X) & \lmt & \psi(\hbar X)H(u,\widetilde{\alpha}(\hbar,X))
\ea$$
is well-defined, smooth and compactly supported. Therefore, I can apply the operator
$D$ to it and get a compactly supported smooth function on $N\times\R\times M$. This
proves that $\widetilde{D}$ is in $\op{\widetilde{G}}{\widetilde{s}}{\R\times M}$.
\end{proof}

Notice that in the above Proof, the function $\psi$ is used to make sense
of expressions of the type $\psi(\hbar X)H(u,\widetilde{\alpha}(\hbar,X))$ even when
$(\hbar,X)$ is not in $\widetilde{W}$, the domain of definition of $\widetilde{\alpha}$.

Let us see what happens to the product of two operators
constructed as in the previous Proposition at $\hbar=0$.
\begin{lemm}\label{lemm:definition-of-tilde}
Let $D_1$ and $D_2$ be in $\op{A}{\tau}{M}$ such that $\widetilde{D_1}$ and
$\widetilde{D_2}$ are in $\mathfrak{Op}(G)$. Let $H$ be a compactly supported
smooth function on $\widetilde{G}$ and $q$ be a point of $M$. If $H_0$ denotes
the restriction of $H$ to $A\subset\widetilde{G}$ then
$$(\widetilde{D_1}\star\widetilde{D_2})_{(0,q)}(H)=(D_1\star D_2)_q(H_0).$$
In particular, if $f_{i=1,2}$ are functions on
$A^*$ such that $\mathfrak{F}(D_{f_1})$ and $\mathfrak{F}(D_{f_2})$,
respectively $\widetilde{\mathfrak{F}(D_{f_1})}$ and $\widetilde{\mathfrak{F}(D_{f_2})}$,
are in $\mathfrak{Op}(A)$, respectively $\mathfrak{Op}(\widetilde{G})$, then
$$(\widetilde{\mathfrak{F}(D_{f_1})}\star\widetilde{\mathfrak{F}(D_{f_2})})_{(0,q)}(H)
=\mathfrak{F}(D_{f_1 f_2})_q(H_0).$$
\end{lemm}
\begin{proof}
The first claim is true because
$$\ba{rcl}
(\widetilde{D_1}\star\widetilde{D_2})_{(0,q)}(H) & = & D_{2,q}(Y\lmt D_{1,q}(X\lmt H(0,X+Y)))\\
 & = & (D_1\star D_2)_q(H_0).
\ea$$
The second claim is true because
$$\ba{cl}
 & (\widetilde{\mathfrak{F}(D_{f_1})}\star\widetilde{\mathfrak{F}(D_{f_2})})_{(0,q)}(H)\\
= & (\mathfrak{F}(D_{f_1})\star\mathfrak{F}(D_{f_2}))_q(H_0)\\
= & \int_{A^*_q}\d\mu(X)f_2(X)\int_{A_q}\d\mu(\xi)e^{-\langle\xi,X\rangle}
\int_{A^*_q}\d\mu(Y)f_1(Y)\int_{A_q}\d\mu(\zeta)e^{-\langle\zeta,Y\rangle}H_0(\xi+\zeta)\\
= & \int_{A^*_q}\d\mu(X)f_2(X)\int_{A_q}\d\mu(\xi)e^{-\langle\xi,X\rangle}
\int_{A^*_q}\d\mu(Y)f_1(Y)\int_{A_q}\d\mu(\zeta)e^{-\langle\zeta-\xi,Y\rangle}H_0(\zeta)\\
= & \int_{A^*_q}\d\mu(X)f_2(X)\mathfrak{F}(\mathfrak{F}\inv(f_1\mathfrak{F}(H_0)))(X)\\
= & \int_{A^*_q}\d\mu(X)f_2(X)f_1(X)\mathfrak{F}(H_0)(X)\\
= & \mathfrak{F}(D_{f_1 f_2})_q(H_0).
\ea$$
\end{proof}

Of course, not every element $D$ of $\op{A}{\tau}{M}$ gives an element of
$\mathfrak{Op}(\widetilde{G})$ and an important problem is to be able
to determine when does $\widetilde{D}$ lie in
$\mathfrak{Op}(\widetilde{G})$? More precisely, for $f$ in $C^\infty(A^*)$, I want
to know when does $\widetilde{\mathfrak{F}(D_f)}$ lie in $\mathfrak{Op}(\widetilde{G})$?

Definition~\ref{defi:When-Is-f-Acceptable-ForQuantisation} gives an answer to this question.

Let $f$ be a smooth function on $A^*$ such that for any $q$ in $M$ and any compactly supported
smooth function $h$ on $A_q$, the product of the restriction $f_q$ of $f$ to $A^*_q$ by
the Fourier transform $\mathfrak{F}(h)$ is again the Fourier transform of a compactly
supported smooth function denoted by $m_f(q)h$
\benn
f_q\mathfrak{F}(h)=\mathfrak{F}(m_f(q)h).
\eenn
For $N$ a smooth manifold and $\theta:N\lra M$ a smooth map, denote by $\Theta$ the induced
bundle morphism $\theta^*A\lra A$.
\begin{defi}\label{defi:When-Is-f-Acceptable-ForQuantisation}
A smooth function $\widetilde{H}$ on $\theta^*A$ is said
to be {\bf sufficiently compact} if
\begin{enumerate}
\item[(i)] it is in $C_{vc}^\infty(\theta^*A)$, the set of vertically
compactly supported smooth functions; and
\item[(ii)] for any subset $K$ of $A$ which is compact
modulo $M$ (that is $K$ is closed and $K\backslash M$ has compact closure)\footnote{
Any compact set is compact modulo $M$; but there might be other compacts modulo $M$:
$M\subset A$ itself is compact modulo $M$ even if it is not compact}, the set
\benn
\big(\supp(\widetilde{H})+\Theta\inv K\big)\cap N
\eenn
is relatively compact. This requirement says that `small vertical perturbations of the support of
$\widetilde{H}$ meet $N$ in a compact set'.
\end{enumerate}
The set of sufficiently compact functions on $\theta^*A$ is denoted by $C_{sc}^\infty(\theta^*A)$.

A smooth function $f$ on $A^*$ is said to be {\bf acceptable for quantisation} if
\begin{enumerate}
\item $f$ has polynomial controlled growth,
\item $m_f$ preserves $C_c^\infty(A)$,\label{item2}
\item $m_f$ preserves $C^\infty_{sc}(\theta^*A)$ for all manifolds $N$ and smooth functions
$\theta:N\lra M$,\label{item3}
\end{enumerate}
where the action of $m_f$ in $(\ref{item2})$ and $(\ref{item3})$ is defined fibrewise.
The set of smooth functions acceptable for quantisation is denoted by $\mathfrak{Q}(A^*)$.
\end{defi}

One reason to state the rather technical above Definition is:
\begin{prop}\label{prop:fInQ(Astar)impliesD_fIsInOp}
If $f$ is acceptable for quantisation then $\widetilde{\mathfrak{F}(D_f)}$ lies in $\mathfrak{Op}
(\widetilde{G})$.
\end{prop}
\begin{proof}
Let $f$ be acceptable for quantisation. Let $N$ be a manifold and $\widetilde{H}$ a compactly
supported smooth function on $N\times\widetilde{G}$. I need to prove that the function
\benn
N\times\widetilde{G} & \lra & \C\\
(u,\hbar,z) & \lmt & \widetilde{\mathfrak{F}(D_f)}_{\hbar,t(z)}\Big((\hbar,x)\lmt
\widetilde{H}(u,(\hbar,x)(\hbar,z))\Big)\\
 & & =\int_{A^*_{t(z)}}\d\mu(X)f(X)\int_{A_{t(z)}}\d\mu(\xi)e^{-i\langle\xi,X\rangle}
\psi(\hbar\xi)\widetilde{H}(u,\widetilde{\alpha}(\hbar,\xi)(\hbar,z)),
\eenn
is compactly supported.

Define $\theta:N\times\widetilde{G}\lra M$ by
$$\theta(u,\hbar,z)=t(z).$$
Let $F$ be the function
\benn
\theta^*A & \lra & \C\\
(u,\hbar,z,\xi) & \lmt & \psi(\hbar\xi)\widetilde{H}(u,\widetilde{\alpha}(\hbar,\xi)(\hbar,z)).
\eenn
I need to prove that $m_f F(u,\hbar,z,0_{t(z)})$ is compactly supported in $(u,\hbar,z)$. It will be
enough to prove that $F$ is sufficiently compact on $\theta^*A$.

Fix $(u,\hbar,z)$ in $N\times\widetilde{G}$. Because $\widetilde{H}$ is compactly supported
and because multiplication on the right in a groupoid is a diffeomorphism between two fibres
of the source map, the map
\benn
N\times\widetilde{s}\inv(\hbar,t(z)) & \lra & \C\\
(u,\hbar,x) & \lmt & \widetilde{H}(u,(\hbar,x)(\hbar,z))
\eenn
is compactly supported. The function $\widetilde{\psi}\circ\widetilde{\alpha}\inv$ is
defined on an open subset of $\widetilde{s}\inv(\hbar,t(z))$ and can be extended by
zero to a smooth function on the whole of $\widetilde{s}\inv(\hbar,t(z))$. Its
product with $\widetilde{H}(u,(\hbar,x)(\hbar,z))$ is compactly supported. This product
composed with $\widetilde{\alpha}$ is a compactly supported function on
$\theta^*A_{(u,\hbar,t(z))}$. This proves that $F$ has vertical compact support.

Let $K$ be a compact modulo $M$ in $A$. I am interested in
\benn
& &(\supp F+\Theta\inv K)\cap N\times\widetilde{G}\\
& \subset & \{(u,\hbar,z,0_{t(z)})\mid
\exists \xi\in A_{t(z)},\,\hbar\xi\in\supp\psi,\,(u,\widetilde{\alpha}(\hbar,\xi)(\hbar,z))
\in\supp \widetilde{H},\,-\xi\in K\}.
\eenn
Let $(u_j,\hbar_j,z_j,0)$ be a sequence in the set on the right hand side of the above
inclusion. For each $j$, choose an element
$\xi_j$ of $-K$ such that $(u_j,\widetilde{\alpha}(\hbar_j,\xi_j)(\hbar_j,z_j))$ is in
the support of $\widetilde{H}$. This sequence satisfies
\begin{enumerate}
\item $-\xi_j$ is a sequence in $K$,
\item $\hbar_j\xi_j$ is a sequence in the support of $\psi$,
\item $(u_j,\widetilde{\alpha}(\hbar_j,\xi_j)(\hbar_j,z_j))$ is a sequence in the
support of $\widetilde{H}$.
\end{enumerate}
Because the support of $\widetilde{H}$ is compact, we can find a subsequence such that
$u_{j_k}$, $\hbar_{j_k}$ and $\widetilde{\alpha}(\hbar_{j_k},\xi_{j_k})(\hbar_{j_k},z_{j_k})$
converge. Also, since $K$ is compact modulo $M$, we can extract a subsequence
such that either $\xi_{j_k}$ converges or $\xi_{j_k}$ lies in $M$. In the former case,
$\hbar_{j_k}\xi_{j_k}$ converges in $\supp\psi$ and $\widetilde{\alpha}(\hbar_{j_k},\xi_{j_k})$
admits a limit, therefore $(\hbar_{j_k},z_{j_k})$ converges. In the latter case,
$\widetilde{\alpha}(\hbar_{j_k},\xi_{j_k})(\hbar_{j_k},z_{j_k})=(\hbar_{j_k},z_{j_k})$ converges
as well. It follows that
$(\supp F+\Theta\inv K)\cap N\times\widetilde{G}$ is compact.
\end{proof}

\begin{prop}\label{prop:a-big-set-of-functions-acceptable-for-quantisation}
If $f$ is either
\begin{itemize}
\item the Fourier transform of a compactly supported smooth function on $A$,
\item polynomial in the fibres, that is $f$ is a smooth section of
$\bigoplus_k S^k A$,
\item a compactly supported character, i.e. it is of the type
$X_q\lmt e^{i\langle\ell(q),X\rangle}$, where $\ell$ is a compactly
supported smooth section of $A$,
\end{itemize}
then $f$ is in $\mathfrak{Q}(A^*)$.
\end{prop}
This Proposition shows that $\mathfrak{Q}(A^*)$ contains indeed many interesting functions.
\begin{proof}
Let $f$ be the Fourier transform of a compactly supported smooth function $g$. For each $q$ in
$M$, choose a local chart together with a
trivialisation of $A$ and $A^*$ over it. Writing things in these local
chart and local trivialisation, to prove that $f$ has polynomial controlled growth is a simple
matter of differentiating under the integral sign in the definition of the Fourier transform.

If $f$ is either polynomial in the fibres or a compactly supported character, it is even
more immediate to prove that $f$ has polynomial controlled growth.

Fix a compactly supported smooth function $H$ on $A$, a manifold $N$, a smooth map
$\theta:N\lra M$ and a sufficiently compact smooth function
$\widetilde{H}\in C_{sc}^\infty(\theta^*A)$.

Firstly, assume that $f$ is the Fourier transform of a compactly supported smooth function
$g$ on $A$. The support of $m_f H$ is included in the sum of $\supp g$ and the support
of $H$, hence it is compactly supported. In the same way, the support of $m_f\widetilde{H}$
is included in the sum of the support of $\widetilde{H}$  and $\Theta\inv\supp g$. It
easily follows that $m_f\widetilde{H}$ is again sufficiently compact.

Secondly, assume that $f$ is a smooth section of $\bigoplus_k S^k A$. For such a function,
the operator $m_f$ is given by a differential operator $\partial_f$ and
$\partial_f H$ has support included in the support
of $H$, therefore $\partial_f H$ is in $C_c^\infty(A^*)$. In the same way, $m_f=\partial_f$
preserves $C_{sc}^\infty(\theta^*A)$. Hence $f$ is acceptable for quantisation.

Finally, assume that $f$ is of the type $f(X_q)=e^{i\langle\ell(q),X\rangle}$.
For $\xi$ in $A_q$, $m_f H(\xi)=H(\xi+\ell(q))$. Hence
$m_f$ translates the support of $H$ by $\ell$ on each fibre of $A\lra M$.
Therefore $m_f H$ has also compact support. For the same reason, $m_f\widetilde{H}$ is
also still vertically compactly supported. The support of $m_f\widetilde{H}$ is equal
to $\supp\widetilde{H}+\Theta\inv\im(-\ell)$. Since $l$ is compactly supported, its image
is compact modulo $M$ and $m_f\widetilde{H}$ is again sufficiently compact.
\end{proof}

\begin{theo}\label{theo:general-quantisation}
Let $\grp{G}{M}$ be a groupoid with Lie algebroid $\tau:A\lra M$. Define a
quantisation map
\benn
\mathcal{Q}:C^\infty_{pg}(A^*) & \lra & \op{\widetilde{G}}{\widetilde{s}}{\R\times M}\\
f & \lmt & \widetilde{\mathfrak{F}(D_f)}.
\eenn
This defines a quantisation of the Poisson manifold $A^*$ in the sense that
$\mathcal{Q}$ sends the set of functions acceptable for quantisation $\mathfrak{Q}(A^*)$
into $\mathfrak{Op}(\widetilde{G})$; moreover, if $f$ and $g$ are two functions acceptable for
quantisation then
\benn
\mathcal{Q}(f)\star\mathcal{Q}(g)_{(0,q)}=\mathcal{Q}(f g)_{(0,q)},
\eenn
and the operator $D=\frac{1}{i\hbar}[\mathcal{Q}(f),\mathcal{Q}(g)]$ is in
$\op{\widetilde{G}}{\widetilde{s}}{\R\times M}$ with
\benn
D_{(0,q)}=\mathcal{Q}(\{f,g\})_{(0,q)},
\eenn
for every $q$ in $M$.
\end{theo}

Notice that along a non zero $\hbar$, $\mathcal{Q}(f)$ restricts to an operator
\benn
\mathcal{Q}(f)_\hbar:C^\infty_c(G)\lra C^\infty_c(G),
\eenn
while for $\hbar=0$ it restricts to an operator
$C^\infty_c(A)\lra C^\infty_c(A)$ which is the Fourier transform of the operator
multiplication by $f$ on $C^\infty_c(A^*)$.

Theorem~\ref{theo:general-quantisation} is a consequence of Lemma~\ref{lemm:definition-of-tilde}
and Corollary~\ref{coro:1-over-hbar-times-the-commutator}. The proof of this
Corollary will take up the whole of next Section.

In Section~\ref{sect:quantisation-of-Rn}, I show how by applying Theorem~\ref{theo:general-quantisation}
one recovers the quantisation of the position and momentum operators used by physicists. In
Section~\ref{sect:the-torus}, I will show how to use it to recover the usual strict
deformation quantisation of the torus with constant Poisson structure.


\section{Computation in local coordinates}\label{sect:Computation-local-coordinates}
Let $m=\dim M$, $U$ be an open subset of $\R^m$ and $\varphi$ a diffeomorphism between
$U$ and an open subset of $M$
\benn
\varphi:U\lra \varphi(U)\subset M.
\eenn
Let
\benn
U\times\R^n & \lra & A\mid_U\\
(u,\xi) & \lmt & \gamma(u,\xi)
\eenn
be a trivialisation above $\varphi(U)$ of the vector bundle $A\lra M$, read in the local
chart $(U,\varphi)$. I identify $\R^n$ with its dual using the usual
euclidean structure of $\R^n$. Therefore $\gamma$ also defines a trivialisation $\delta$
of the restriction of $A^*\lra M$ to $\varphi(U)$. This trivialisation is characterised by
\benn
\langle\delta(u,X),\gamma(u,\xi)\rangle=\langle X,\xi\rangle,
\eenn
where $\langle\,,\rangle$ denotes both the pairing between $A$ and $A^*$, and the euclidean
product on $\R^n$.

Choose an open neighbourhood $V_U$ of $U\times\{0\}$ in $\R^m\times\R^n$,
such that $\alpha$ is defined on
$\gamma(V_U)$. I can define a local chart for $G$
\benn
\theta: V_U & \lra & G\\
(u,v) & \lmt & \alpha\circ\gamma(u,v).
\eenn
Let $V_U^\prime$ be an open neighbourhood of $U\times\{0\}$ in $V_U$ verifying
\begin{enumerate}
\item for each $(u,v)$ in $V_U^\prime$, there exists $u$ in $U$ such that
$t\circ\theta(u,v)=\varphi(u)$,
\item for each $(u_1,v_1)$ and $(u_2,v_2)$ in $V_U^\prime$ with $t\circ\theta(u_2,v_2)
=\varphi(u_1)$, the product $\theta(u_1,v_1)\theta(u_2,v_2)$ is in $\theta(V_U)$.
\end{enumerate}
Let $\sigma:V_U^\prime\lra U$ be given by
\benn
\sigma(u,v)=\varphi\inv\circ t\circ\theta(u,v),
\eenn
(notice that $\varphi\inv\circ s\circ\theta(u)=u$). Let
\benn
V_U^\prime\,_{\mbox{pr}_U}\!\!\times_\sigma V_U^\prime=\{(u,v_1,v_2)\mid\,(u,v_2)\in V_U^\prime,\,
(\sigma(u,v_2),v_1)\in V_U^\prime\}
\eenn
and define
\benn
p:V_U^\prime\,_{\mbox{pr}_U}\!\!\times_\sigma V_U^\prime & \lra & V_U\\
(u,v_1,v_2) & \lmt & \theta\inv(\theta(\sigma(u,v_2),v_1)\theta(u,v_2)).
\eenn
Ramazan~\cite[Proposition $2.2.5$]{Ra} proved
\benn
p(u,v_1,v_2) & = & (u,v_1+v_2+B(u,v_1,v_2)+\mbox{O}_3(u,v_1,v_2)),\\
\theta(u,v)\inv & = & \theta(\sigma(u,v),-v+B(u,v,v)+\mbox{O}_3(u,v)),
\eenn
where $B(u,v_1,v_2)$ is bilinear in $(v_1,v_2)$ and $\mbox{O}_3(u,v_1,v_2)$,
respectively $\mbox{O}_3(u,v)$, is
of degree of homogeneity at least $3$ in $v_1$ and $v_2$, respectively $v$.

For $\xi$ in $\R^n$
\benn
\frac{\partial\theta}{\partial v}(u,v)\xi=\d_{\gamma(u,v)}\alpha\circ\frac{\partial\gamma}
{\partial v}(u,v)\xi=\d_{\gamma(u,v)}\alpha\circ\gamma(u,\xi),
\eenn
because $\gamma$ is linear in $v$. In particular
\benn
\frac{\partial\theta}{\partial v}(u,0)\xi=\gamma(u,\xi),
\eenn
because $\d\alpha$ is the identity along $M$. Moreover $\varphi(u)=\theta(u,0)$.

The map
\benn
\widetilde{\gamma}:\R\times U\times\R^n & \lra & \widetilde{A}\mid_{\R\times U}\\
(\hbar,u,\xi) & \lmt & (\hbar,\gamma(u,\xi))
\eenn
gives a local trivialisation of the Lie algebroid $\widetilde{A}\lra \R\times M$ over
$\R\times U$. Let
\benn
\widetilde{V_U}=\{(\hbar,u,v)\mid\,(u,\hbar v)\in V_U\}
\eenn
and
\benn
\widetilde{V_U^\prime}=\{(\hbar,u,v)\mid\,(u,\hbar v)\in V_U^\prime\}.
\eenn
I obtain local coordinates on $\widetilde{G}$ by taking $\widetilde{\theta}=
\widetilde{\alpha}\circ\widetilde{\gamma}$
\benn
\widetilde{\theta}:\widetilde{V_U} & \lra & \widetilde{G}\\
(\hbar,u,v) & \lmt & \widetilde{\alpha}(\hbar,\gamma(u,v)).
\eenn

Let $q$ be in $M$ and $\xi$ be in $A_q$. Assume that $\xi$ is in the domain
of $\alpha$. The map
\benn
\Th_{\xi}:A_q & \lra & A_{t\circ\alpha(\xi)}\\
\zeta & \lmt & \left.\frac{\d}{\d r}\right|_{r=0}\alpha(\xi+r\zeta)\alpha(\xi)\inv
\eenn
defines an isomorphism between $A_q$ and $A_{t\circ\alpha(\xi)}$.

\begin{lemm}\label{lemm:crucialdifhz}
Let $u$ be in $U$. Let $\hbar\in\R$ and $\zeta$, $\xi$ in $\R^n$
\benn
\difhz\widetilde{\alpha}(\hbar,\Th_{\hbar\gamma(u,\xi)}\circ\gamma(u,\zeta))\widetilde{\alpha}
(\hbar,\gamma(u,\xi))=\d_{(0,u,\zeta+\xi)}\widetilde{\theta}(1,0,0).
\eenn
\end{lemm}
\begin{proof}
I first compute

\benn
 & & \widetilde{\alpha}(\hbar,\left.\frac{\d}{\d r}\right|_{r=0}
 \alpha\circ\gamma(u,\hbar\xi+r\zeta)
 \alpha\circ\gamma(u,\hbar\xi)\inv)\widetilde{\alpha}(\hbar,\gamma(u,\xi))\\
 & = & \widetilde{\alpha}(\hbar,\left.\frac{\d}{\d r}\right|_{r=0}
 \theta(u,\hbar\xi+r\zeta)
 \theta(u,\hbar\xi)\inv)\widetilde{\theta}(\hbar,u,\xi)\\
 & = & \widetilde{\alpha}\Big(\hbar,\left.\frac{\d}{\d r}\right|_{r=0}
 \theta(u,\hbar\xi+r\zeta)
 \theta(\sigma(u,\hbar\xi),-\hbar\xi+\mbox{O}(\hbar^2))\Big)
 \widetilde{\theta}(\hbar,u,\xi)\\
 & = & \widetilde{\alpha}\Big(\hbar,\left.\frac{\d}{\d r}\right|_{r=0}
 \theta(\sigma(u,\hbar\xi),r\zeta+\mbox{O}(\hbar^2)+B(\sigma(u,\hbar\xi),\hbar\xi+r\zeta,
 -\hbar\xi+\mbox{O}(\hbar^2))\\
 & & +\mbox{O}_3(\sigma(u,\hbar\xi),\hbar\xi+r\zeta,
 -\hbar\xi+\mbox{O}(\hbar^2)))\Big)\widetilde{\theta}(\hbar,u,\xi)\\
 & = & \widetilde{\theta}\circ\widetilde{\gamma}\inv\Big(\hbar,\frac{\partial\theta}{\partial v}(\sigma(u,\hbar\xi),0)(
 \zeta+B(\sigma(u,\hbar\xi),\zeta,
 -\hbar\xi+\mbox{O}(\hbar^2))+\mbox{O}(\hbar^2))\Big)\widetilde{\theta}(\hbar,u,\xi)
 \footnotemark\\
 & = & \widetilde{\theta}\Big(\hbar,\sigma(u,\hbar\xi),
 \zeta-\hbar B(\sigma(u,\hbar\xi),\zeta,\xi)+\mbox{O}(\hbar^2)\Big)\widetilde{\theta}(\hbar,u,\xi)\\
 & = & (\hbar,\theta(\sigma(u,\hbar\xi),
 \hbar\zeta-\hbar^2 B(\sigma(u,\hbar\xi),\zeta,\xi)+\mbox{O}(\hbar^3))\theta(u,\hbar\xi))
 \footnotemark\\
 & = & (\hbar,\theta(u,\hbar\zeta-\hbar^2 B(\sigma(u,\hbar\xi),\zeta,\xi)
 +\mbox{O}(\hbar^3)+\hbar\xi+\\
 & & +B(u,\hbar\zeta-\hbar^2B(\sigma(u,\hbar\xi),\zeta,\xi)
 +\mbox{O}(\hbar^3),\hbar\xi)+\mbox{O}(\hbar^3)))\\
 & = & (\hbar,\theta(u,\hbar\zeta+\hbar\xi-\hbar^2 B(\sigma(u,\hbar\xi),\zeta,\xi)
 +B(u,\hbar\zeta,\hbar\xi)+\mbox{O}(\hbar^3)))\\
 & = & \widetilde{\theta}(\hbar,u,\zeta+\xi-\hbar B(\sigma(u,\hbar\xi),\zeta,\xi)
 +\hbar B(u,\zeta,\xi)+\mbox{O}(\hbar^2))).
\eenn
\addtocounter{footnote}{-1}\footnotetext{At $r=0$, we have $\theta(u,\hbar\xi+r\zeta)\theta(u,\hbar\xi)\inv=t\circ\theta(u,\xi)$, thus
 the differential of $\theta(u,\hbar\xi+r\zeta)\theta(u,\hbar\xi)\inv$ at $r=0$ is of the
 type $\frac{\partial\theta}{\partial v}(\sigma(u,\hbar\xi),0)\phi$ for a certain vector
 $\phi$.}
\stepcounter{footnote}
\footnotetext{This is true only for $\hbar\neq0$, nevertheless the final result of
the computation is trivially true for $\hbar=0$.}
The Lemma follows by differentiation with respect to $\hbar$ at $0$.
\end{proof}

\begin{lemm}\label{lemm:diff-de-N-fg}
Let $f$ and $g$ be in $\mathfrak{Q}(A^*)$. Let
$q$ be a point in $M$ and $H$ a compactly supported smooth function on $\widetilde{G}$. Let
\benn
N_{f g}(\hbar)=\widetilde{\mathfrak{F}(D_f)}\star\widetilde{\mathfrak{F}(D_g)}_{\hbar,q}(H),
\eenn
then
\benn
\hspace*{-2cm}\frac{\d N_{f g}}{\d\hbar}(0)=
\mathfrak{F}(D_{f^\prime g})_{q}(H_0)
+\mathfrak{F}(D_{f g})_{q}(\left.\frac{\partial H}{\partial\hbar}\right|_{\hbar=0}),
\eenn
where
$$f^\prime(Y)=\difhz f(Y\circ\Th_{\hbar\xi}\inv),$$
and $H_0=H\mid_{\hbar=0}$.
\end{lemm}
The term $\left.\frac{\partial H}{\partial\hbar}\right|_{\hbar=0}$ is defined by first
pulling back $H$ to a neighbourhood of $\{0\}\times A$ in $\R\times A$ via $\widetilde{\alpha}$,
then differentiating with
respect to $\hbar$ and finally pushing forward the result via $\widetilde{\alpha}$
again. This definition is actually independent of the choice of $\alpha$.
\begin{proof}
We have
\benn
N_{f g}(\hbar) & = & \int_{A^*_q}\d\mu(X)g(X)\int_{A_q}\d\mu(\xi)e^{-i\langle X,\xi\rangle}
\int_{A^*_{t\circ\alpha(\hbar\xi)}}\d\mu(Y)f(Y)\\
 & & \int_{A_{t\circ\alpha(\hbar\xi)}}\d\mu(\zeta)
e^{-i\langle Y,\zeta\rangle}\psi(\hbar\xi)\psi(\hbar\zeta)H(\widetilde{\alpha}(\hbar,\zeta)\widetilde{\alpha}(\hbar,\xi))
\eenn
The following change of variables,
\begin{itemize}
\item replace $Y$ by $Y\circ\Th_{\hbar\xi}\inv$ with $Y\in A_q^*$,
\item replace $\zeta$ by $\Th_{\hbar\xi}(\zeta)$ with $\zeta\in A_q$,
\end{itemize}
gives
\benn
N_{f g}(\hbar) & = & \int_{A^*_q}\d\mu(X)g(X)\int_{A_q}\d\mu(\xi)e^{-i\langle X,\xi\rangle}
\int_{A^*_q}\d\mu(Y)f(Y\circ\Th_{\hbar\xi}\inv)\\
 & & \int_{A_q}\d\mu(\zeta)
e^{-i\langle Y,\zeta\rangle}\psi(\hbar\xi)\psi(\hbar\Th_{\hbar\xi}(\zeta))H(\widetilde{\alpha}(\hbar,\Th_{\hbar\xi}(\zeta))\widetilde{\alpha}(\hbar,\xi)).
\eenn
The changes of variables require to introduce the terms $\det\Th_{\hbar\xi}$ and
$\det\Th_{\hbar\xi}\inv$ in the above integral; but these two terms cancel each other
since their product is $1$.

Leaving out the justification for it for later, I differentiate the above expression
under the integral signs. Since $\psi$ is constant and equal to $1$ in a
neighbourhood of $M$, it follows that
\begin{itemize}
\item $\psi(0_q)=1$ and
\item $\difhz\psi(\hbar\xi)=\difhz\psi(\hbar\Th_{\hbar\xi}(\zeta))=0$.
\end{itemize}
Because of Lemma~\ref{lemm:crucialdifhz} and by definition of
$\frac{\partial H}{\partial\hbar}$
\benn
\difhz H(\widetilde{\alpha}(\hbar,\Th_{\hbar\xi}(\zeta))\widetilde{\alpha}(\hbar,\xi))=
\frac{\partial H}{\partial\hbar}(0,\zeta+\xi).
\eenn
The Lemma follows since
\benn
 & & \int_{A^*_q}\d\mu(X)g(X)\int_{A_q}\d\mu(\xi)e^{-i\langle X,\xi\rangle}
 \int_{A^*_q}\d\mu(Y)f(Y)\int_{A_q}\d\mu(\zeta)
 e^{-i\langle Y,\zeta\rangle}\frac{\partial H}{\partial\hbar}(0,\zeta+\xi)\\
 & = & \mathfrak{F}(D_f)\star\mathfrak{F}(D_g)_q
 (\left.\frac{\partial H}{\partial\hbar}\right|_{\hbar=0})\\
 & = & \mathfrak{F}(D_{f g})_q(\left.\frac{\partial H}{\partial\hbar}\right|_{\hbar=0}),
\eenn
where the last line is true by Lemma~\ref{lemm:definition-of-tilde}.

There now remains to justify differentiation below the integral signs in
\benn
& & \int_{A^*_q}\d\mu(X)g(X)\int_{A_q}\d\mu(\xi)e^{-i\langle X,\xi\rangle}
\int_{A^*_q}\d\mu(Y)f(Y\circ\Th_{\hbar\xi}\inv)\\
 & & \int_{A_q}\d\mu(\zeta)
e^{-i\langle Y,\zeta\rangle}\psi(\hbar\xi)\psi(\hbar\Th_{\hbar\xi}(\zeta))
H(\widetilde{\alpha}(\hbar,\Th_{\hbar\xi}(\zeta))\widetilde{\alpha}(\hbar,\xi)).
\eenn
Let $\theta$ be the map
\benn
\widetilde{V_U} & \lra & M\\
(\hbar,\xi) & \lmt & t\circ\alpha(\hbar\xi).
\eenn
Define a function $\widetilde{H}$ on $\theta^*A$ by
\benn
\widetilde{H}(\hbar,\xi,\zeta)=\psi(\hbar\xi)\psi(\hbar\zeta)H(\widetilde{\alpha}(\hbar,\zeta)
\widetilde{\alpha}(\hbar,\xi)).
\eenn
Let
\benn
S_1(\hbar,\xi,Y)& = & \int_{A_q}\d\mu(\zeta)e^{-i\langle Y,\zeta\rangle}\widetilde{H}(\hbar,\xi,\zeta),\\
S_2(\hbar,\xi) & = & \int_{A^*_q}\d\mu(Y) f(Y\circ\Th\inv_{\hbar\xi})S_1(\hbar,\xi,Y),\\
S_3(\hbar,X) & = & \int_{A_q}\d\mu(\xi) e^{-i\langle X,\xi\rangle}S_2(\hbar,\xi),
\eenn
where $S_1$ is defined on $\theta^*A^*$, $S_2$ is defined on $\widetilde{V_U}$ and extended
by zero to $\R\times A$ and $S_3$ is defined on $\R\times A^*$.

I claim that $\widetilde{H}$ is in $C^\infty_{sc}(\theta^*A)$. The proof of this claim is similar
to that for $F$ in the proof of Proposition~\ref{prop:fInQ(Astar)impliesD_fIsInOp}
and will not be reproduced here. It follows that
\benn
S_2(\hbar,\xi) & = & m_f\widetilde{H}(\hbar,\xi,0)
\eenn
is compactly supported in $(\hbar,\xi)$.

For $\xi$ fixed, the function $\widetilde{H}$ is compactly supported in $(\hbar,\zeta)$, hence
derivation below the integral sign in $S_1$ is possible.

Since $f$ has polynomial controlled growth, for $\xi$ fixed, there exists a positive
constant $C$, an $\epsilon>0$ and an integer $m$ such that
\benn
\left|\frac{\partial}{\partial\hbar}\big(f(Y\circ\Th\inv_{\hbar\xi})S_1(\hbar,\xi,Y)\big)\right|\leq
C(1+\|Y\|^2)^m(|S_1(\hbar,\xi,Y)|+|\frac{\partial}{\partial\hbar}S_1(\hbar,\xi,Y)|).
\eenn
By Lemma~\ref{lemm:Fourier-transform-of-compactly-supp-is-Schwartz},
both terms on the right hand side of the above inequality are bounded by a smooth
$L^1$ function independent of $\hbar$. Differentiation below the integral sign in
$S_2$ is therefore possible.

Since $S_2$ is compactly supported, differentiation below the integral sign in $S_3$ is
possible.

To finish, since $g$ has polynomial controlled growth and by
Lemma~\ref{lemm:Fourier-transform-of-compactly-supp-is-Schwartz}, differentiation
below the integral sign in $N_{f g}$ is possible.
\end{proof}
\begin{lemm}\label{lemm:diff-of-Thbarxi-composed-with-gamma}
Let $\xi$ and $\zeta$ be in $\R^n$, then
\benn
\difhz\Th_{\gamma(u,\hbar\xi)}\circ\gamma(u,\zeta)=\frac{\partial\gamma}{\partial u}
(u,\zeta)\circ\frac{\partial\sigma}{\partial v}(u,0)\xi-\gamma(u,B(u,\zeta,\xi)).
\eenn
\end{lemm}
\begin{proof}
First, I compute
\benn
 & & \Th_{\gamma(u,\hbar\xi)}\circ\gamma(u,\zeta)\\
 & = & \left.\frac{\d}{\d r}\right|_{r=0}\theta(u,\hbar\xi+r\zeta)\theta(\sigma(u,\hbar\xi),
 -\hbar\xi+\mbox{O}(\hbar^2))\\
 & = & \left.\frac{\d}{\d r}\right|_{r=0}\theta(\sigma(u,\hbar\xi),r\zeta+
 B(\sigma(u,\hbar\xi),\hbar\xi+r\zeta,-\hbar\xi+\mbox{O}(\hbar^2))+\\
 & & +\mbox{O}_3
 (\sigma(u,\hbar\xi),\hbar\xi+r\zeta,-\hbar\xi+\mbox{O}(\hbar^2)))\\
 & = & \frac{\d\theta}{\d v}(\sigma(u,\hbar\xi),0)(\zeta-\hbar B(\sigma(u,\hbar\xi),\zeta,\xi)+
 \mbox{O}(\hbar^2))\\
 & = & \gamma(\sigma(u,\hbar\xi),\zeta-\hbar B(\sigma(u,\hbar\xi),\zeta,\xi)+\mbox{O}(\hbar^2)).
\eenn
The result follows by differentiation and because $\gamma$ is linear in the second variable.
\end{proof}
The map $\gamma:U\times\R^n\lra A\mid_U$ is a local trivialisation of $A$. The induced local
trivialisation of $A^*$
\benn
\delta:U\times\R^n & \lra & A^*\mid_U
\eenn
is characterised by
\benn
\langle\delta(u,X),\gamma(u,\zeta)\rangle=\langle X,\zeta\rangle,
\eenn
where $\langle\,,\rangle$ denotes the euclidean product on $\R^n$.
\begin{lemm}\label{lemm:delta-inv-Y-T-kh-inv}
Let $e_1,\ldots,e_n$ be a basis of $\R^n$, then for $(u,\xi)$ in $U\times \R^n$ and $Y$ in
$A^*_{\varphi(u)}$
\benn
\difhz\delta\inv(Y\circ\Th_{\hbar\gamma(u\xi)}\inv)=\Big(\frac{\partial\sigma}{\partial v}(u,0)
\xi,\sum_k Y\big(\gamma(u,B(u,e_k,\xi))\big)e_k\Big).
\eenn
\end{lemm}
\begin{proof}
In the proof of Lemma~\ref{lemm:diff-of-Thbarxi-composed-with-gamma}, I showed that
\benn
\Th_{\gamma(u,\hbar\xi)}\circ\gamma(u,\zeta)=
\gamma(\sigma(u,\hbar\xi),\zeta-\hbar B(\sigma(u,\hbar\xi),\zeta,\xi)+\mbox{O}(\hbar^2)),
\eenn
therefore
\benn
\Th_{\gamma(u,\hbar\xi)}\inv\circ\gamma(\sigma(u,\hbar\xi),\zeta)=
\gamma(u,\zeta+\hbar B(u,\zeta,\xi)+\mbox{O}(\hbar^2)).\footnotemark
\eenn
\footnotetext{Here I used $\sigma(u,\hbar\xi)=u+\mbox{O}(\hbar)$ to simplify the formula.}
The Lemma is proved by using this formula when differentiating
\benn
\delta\inv(Y\circ\Th_{\hbar\gamma(u,\xi)}\inv)=\Big(\sigma(u,\hbar\xi),
\sum_k Y\big(\Th_{\hbar\gamma(u\xi)}\inv\circ\gamma(\sigma(u,\hbar\xi),e_k)\big)e_k\Big).
\eenn
\end{proof}
Write
\benn
B(u,e_k,e_h)=\sum_j B_{k,h}^j e_j.
\eenn
In particular, the $B_{k h}^j$'s depend on $u$.
\begin{coro}\label{coro:diff-de-f}
For $Y=\sum_j Y_j e_j$ and $\xi=\sum_h\xi_h e_h$ in $\R^n$
\benn
\difhz\delta\inv(\delta(u,Y)\circ\Th_{\hbar\gamma(u,\xi)}\inv)=\Big(\frac{\partial\sigma}{\partial v}(u,0)
\xi,\sum_{k,h,j} Y_j B_{k,h}^j\xi_h e_k\Big).
\eenn
Let $f$ be a smooth function on $A^*$. Define $F=f\circ\delta$, then
\benn
\difhz f(\delta(u,Y)\circ\Th_{\hbar\gamma(u,\xi)}\inv)=\sum_{k,h}\xi_k\frac{\partial F}
{\partial u_h}(u,Y)\frac{\partial\sigma_h}{\partial v_k}(u,0)+\sum_{k,h,j}
Y_j\xi_h\frac{\partial F}{\partial Y_k}(u,Y)B^j_{k,h}.
\eenn
\end{coro}
\begin{proof}
The first formula is just Lemma~\ref{lemm:delta-inv-Y-T-kh-inv} written in local
coordinates. The second one is a straightforward computation.
\end{proof}

Let us look at the Poisson bracket on $A^*$ in local coordinates.
\begin{lemm}\label{lemm:poisson-bracket-of-f-and-g}
Let $f$ and $g$ be smooth functions on $A^*$. Define $F=f\circ\delta$ and $G=g\circ\delta$,
smooth functions on $U\times\R^n$. Set
\benn
\{F,G\}=\{f,g\}\circ\delta.
\eenn
Let $(u,Z)$ be in $U\times\R^n$ and denote $\frac{\partial F}{\partial u_j}$
for $\frac{\partial F}{\partial u_j}(u,Z)$. I will use similar notations for
$\frac{\partial G}{\partial u_j}$, $\frac{\partial F}{\partial Z_k}$ and
$\frac{\partial G}{\partial Z_k}$. Then
\benn
\{F,G\}(u,Z)=\sum_{k,h,j}\frac{\partial F}{\partial Z_k}\frac{\partial G}{\partial Z_h}
(B^j_{h k}-B^j_{k h})Z_j\,+\,\sum_{k,h}(
\frac{\partial F}{\partial Z_k}\frac{\partial G}{\partial u_h}-
\frac{\partial F}{\partial u_h}\frac{\partial G}{\partial Z_k})
\frac{\partial \sigma_h}{\partial v_k}(u,0).
\eenn
\end{lemm}
\begin{proof}
This is essentially Equation~$(1.2.6)$ and Proposition~$2.2.6$ in Ramazan~\cite{Ra}
where it is proved that\footnote{The signs here and in Ramazan~\cite{Ra} do not
agree. This is due to different choices in the definition of the map $\alpha$. Essentially,
I have $s\circ\alpha$ constant on the fibres of $A\lra Q$, whereas he has $t\circ\alpha$
constant on the same fibres.}
\benn
[\gamma(u,e_k),\gamma(u,e_h)]=\sum_j (B^j_{k h}-B^j_{h k})\gamma(u,e_j)
\eenn
and, if $u_j^*$ is the $j$-th coordinate map on $U$,
\benn
\rho(e_k)\cdot u_h^*(u)=\frac{\partial\sigma_h}{\partial v_k}(u,0).
\eenn
\end{proof}

\begin{prop}\label{prop:differential-of-the-commutator}
Let $f$ and $g$ be in $\mathfrak{Q}(A^*)$. Let $H$ be a compactly supported smooth function
on $\widetilde{G}$ and $q$ a point in $M$, then
\benn
\difhz\Big(\widetilde{\mathfrak{F}(D_f)}\star\widetilde{\mathfrak{F}(D_g)}
-\widetilde{\mathfrak{F}(D_g)}\star\widetilde{\mathfrak{F}(D_f)}\Big)_{(\hbar,q)}(H)=
i\widetilde{\mathfrak{F}(D_{\{f,g\}})}_{(0,q)}(H).
\eenn
\end{prop}
\begin{proof}
The left hand side of the above equation is equal to
\benn
\difhz(N_{f g}-N_{g f}).
\eenn
The terms
$$\mathfrak{F}(D_{f g})_{q}(\left.\frac{\partial H}{\partial\hbar}\right|_{\hbar=0})$$
appear in both $\difhz N_{f g}$ and $\difhz N_{g f}$ in Lemma~\ref{lemm:diff-de-N-fg};
they will therefore cancel each other when taking the difference.
The other term in $\difhz N_{f g}$, when using Corollary~\ref{coro:diff-de-f},
becomes a sum of terms. These terms can be dealt with by recalling that the Fourier
transform takes the operator `multiplication by a variable' to the operator 
`derivation with respect to this variable'. For example (with some slight abuse of notations)
\benn
 & &\int_{\R^n}\d\mu(X)G(u,X)\int_{\R^n}\d\mu(\xi)e^{-i\langle X,\xi\rangle}\xi_k\int_{\R^n}\d\mu(Y)
\frac{\partial F}{\partial u_h}(u,Y)\frac{\partial\sigma_h}{\partial v_k}(u,0)\\
& & \int_{\R^n}\d\mu(\zeta)e^{-i\langle Y,\zeta\rangle}H_0(\xi+\zeta)\\
& = & \int_{\R^n}\d\mu(X)G(u,X)i\frac{\partial}{\partial X_k}
\int_{\R^n}\d\mu(\xi)e^{-i\langle X,\xi\rangle}\int_{\R^n}\d\mu(Y)
\frac{\partial F}{\partial u_h}(u,Y)\frac{\partial\sigma_h}{\partial v_k}(u,0)\\
& & \int_{\R^n}\d\mu(\zeta)e^{-i\langle Y,\zeta\rangle}H_0(\xi+\zeta)\\
& = & -i\int_{\R^n}\d\mu(X)\frac{\partial G}{\partial X_k}(u,X)
\int_{\R^n}\d\mu(\xi)e^{-i\langle X,\xi\rangle}\int_{\R^n}\d\mu(Y)
\frac{\partial F}{\partial u_h}(u,Y)\frac{\partial\sigma_h}{\partial v_k}(u,0)\\
& & \int_{\R^n}\d\mu(\zeta)e^{-i\langle Y,\zeta\rangle}H_0(\xi+\zeta)\\
& = & -i\int_{\R^n}\d\mu(Z)\frac{\partial G}{\partial Z_k}(u,Z)
\frac{\partial F}{\partial u_h}(u,Z)\frac{\partial\sigma_h}{\partial v_k}(u,0)\int_{\R^n}\d\mu(\xi)
e^{-i\langle Z,\xi\rangle}H_0(\xi),
\eenn
where the last line is true for the same reason that $\mathfrak{F}(D_f)\star
\mathfrak{F}(D_g)=\mathfrak{F}(D_{f g})$ (see Lemma~\ref{lemm:definition-of-tilde}).
A similar computation leads to
\benn
 & &\int_{\R^n}\d\mu(X)G(u,X)\int_{\R^n}\d\mu(\xi)e^{-i\langle X,\xi\rangle}\xi_h\int_{\R^n}\d\mu(Y)
Y_j B_{k h}^j\frac{\partial F}{\partial Y_k}(u,Y)
\int_{\R^n}\d\mu(\zeta)e^{-i\langle Y,\zeta\rangle}H_0(\xi+\zeta)\\
& = & -i\int_{\R^n}\d\mu(Z)\frac{\partial G}{\partial Z_h}(u,Z)
\frac{\partial F}{\partial Z_k}(u,Z)Z_j B^j_{k h}\int_{\R^n}\d\mu(\xi)
e^{-i\langle Z,\xi\rangle}H_0(\xi).
\eenn
These computations together with Lemma~\ref{lemm:poisson-bracket-of-f-and-g} prove
Proposition~\ref{prop:differential-of-the-commutator}.
\end{proof}
As a corollary, I obtain the following:
\begin{coro}\label{coro:1-over-hbar-times-the-commutator}
Let $f$ and $g$ be two smooth functions on $A^*$ acceptable for quantisation. Then
\benn
\frac{1}{i\hbar}[\widetilde{\mathfrak{F}(D_f)},\widetilde{\mathfrak{F}(D_g)}]
\eenn
is a well-defined element of $\op{\widetilde{G}}{\widetilde{s}}{\R\times M}$, which along
$\hbar=0$ is equal to
\benn
\widetilde{\mathfrak{F}(D_{\{f,g\}})}_0.
\eenn
\end{coro}

\section{Quantisation of $\R^{2n}$}\label{sect:quantisation-of-Rn}
In this short section, I will discuss the case of the quantisation of observables on
the phase space of a particle moving in $\R^n$.

Let $M=\R^n$ with its euclidean structure and $A$ be the tangent bundle of $M$, that
is $A=\R^{2n}$. With these conditions, the Lie groupoid $G$ integrating $A$ is the
pair groupoid $M\times M=\R^n\times\R^n$ with source map, tangent map and product
\benn
s(p,q)=q,\,t(p,q)=p,\,(r,p)\cdot(p,q)=(r,q).
\eenn
The euclidean product gives a natural family of measures on the fibres of $A$. I can
take $\alpha$ to be defined on the whole of $A$ by
\benn
\alpha(q,\xi)=(q+\xi,q).
\eenn
Also, I can choose $\psi$ to be equal to the constant function $1$. The space of
morphisms of the tangent groupoid is diffeomorphic to $\R\times\R^n\times\R^n$.
Let $H$ be a compactly supported smooth function on $\R\times\R^n\times\R^n$; it is a
a function of $(\hbar,p,q)$.
Let $f$ be a function acceptable for quantisation. For $\hbar\neq 0$ and $x=(p,q)$ in $G$
\benn
\mathcal{Q}_\hbar(f)(H)(x) & = & \int_{A^*_p}\d X\,f(p,X)\int_{A_p}
e^{-i\langle\xi,X\rangle}H(\widetilde{\alpha}(\hbar,p,\xi)\cdot(\hbar,p,q))\\
 & = & \int_{\R^n}\d X\,f(p,X)\int_{\R^n}
e^{-i\langle\xi,X\rangle}H((\hbar,p+\hbar\xi,p)\cdot(\hbar,p,q))\\
 & = & \int_{\R^n}\d X\,f(p,X)\int_{\R^n}
e^{-i\langle\xi,X\rangle}H(\hbar,p+\hbar\xi,q).
\eenn
If $f(p,X)$ is equal to $X_k$, the $k$-th coordinate of $X$, its quantisation at a given value
of $\hbar$ is
\benn
\mathcal{Q}_\hbar(X_k)=-i\hbar\frac{\partial}{\partial p_k};
\eenn
whereas if $f(p,X)=p_k$ then
\benn
\mathcal{Q}_\hbar(p_k)=p_k,
\eenn
the operator multiplication by $p_k$.


\section{Some general results about Poisson manifolds}\label{sect:general-results-about-poisson}
Given a Poisson manifold $P$, its cotangent bundle is naturally a Lie algebroid $A$. If
this Lie algebroid is integrable to a Lie groupoid then Theorem~\ref{theo:general-quantisation}
gives a quantisation of the Lie-Poisson manifold $A^*$. Since the Poisson structure of $A^*$
is completely determined by the one of $P$, and vice-versa, one might hope to
be able to say something about
a quantisation of $P$. One way of doing so might consist in looking for a surjective
Poisson map $\pi$ between $A^*$ and $P$ and then quantised a function on $P$ by taking the
quantisation of the pulled back function on $A^*$. In $C^*$-algebraic quantisation, such
an idea is bound to fail because if $f$ is a function on $P$, then its pull-back
$\pi^*f$ has little chance of being quantisable. Nevertheless, I will show with an example
that this idea can be made to work when using the quantisation defined in
Theorem~\ref{theo:general-quantisation}.

The aim of this Section is to derive the partial differential equation that a map
$\pi:T P\lra P$ has to satisfy to be Poisson. This PDE is given in
Corollary~\ref{coro:poissonmapforriemanniancompatiblestructure}.

Assume $P$ is a Poisson manifold with Poisson bivector field $\eta$.
Denote by $p:\T P\rightarrow P$ the natural projection.
The Poisson structure of $P$ induces a Lie algebroid structure on the cotangent space
of $P$ with anchor map $\eta:T^*P\lra T P$\footnote{Here, $\eta$ is understood as an
anti-symmetric map
between $T^*P$ and $T P$. I will use different sorts of interpretations of $\eta$, the precise
interpretation depending on the context.}.
Its dual, the tangent space of $P$, inherits a Poisson structure in the
following manner. Let $\alpha$ and $\beta$ be closed $1$-forms on $P$. They
naturally define smooth functions, denoted $\widetilde{\alpha}$ and $\widetilde{\beta}$,
on $T P$ by, for $v$ in $T_x P$,
$$\widetilde{\alpha}(v)=\alpha(v)\mbox{ and }\widetilde{\beta}(v)=\beta(v).$$
Put
\be
\{\widetilde{\alpha},\widetilde{\beta}\}_{\T P}(v) & = & [\alpha,\beta](v), \label{eqn:bracketonTV1}
\ee
where $[\,,]$ is the bracket on the Lie algebroid $\T^*P$.
Let $f$ and $g$ be smooth functions on $P$. Put
\be
\{p^*f,p^*g\}_{\T P} & = & 0. \label{eqn:bracketonTV2}
\ee
Finally, for $v$ in $T_x P$, put
\be
\{\widetilde{\alpha},p^*f\}_{\T P}(v) & = & \eta(\alpha(x))\cdot f. \label{eqn:bracketonTV3}
\ee
The formulae (\ref{eqn:bracketonTV1}), (\ref{eqn:bracketonTV2}) and (\ref{eqn:bracketonTV3})
completely determine the Poisson structure on $\T P$.

Suppose there is a torsion free connection on $T P$
\benn
\nabla:\Gamma(T P)\otimes\Gamma(T P) & \lra & \Gamma(T P)\\
(X,Y) & \lmt & \nabla_X Y.
\eenn
For example, the Levi-Civita
connection of a metric would do. In particular,
for every $v$ in $\T_x P$, there is a splitting of $\T_v(\T P)$ as a direct sum of a horizontal
space $\H_v(\T P)$ and a vertical space $\V_v(\T P)$. Denote by $\H$ and $\V$ the projection
on respectively $\H(\T P)$ and $\V(\T P)$. Both these spaces are isomorphic to $\T_x P$ and the
projection $\H$ is equal to $p_*$ while $\V(\T P)$ is the kernel of $p_*$. The
isomorphism between $\V_v(\T P)$ and $\T_x P$ is given by
\benn
\T_x P & \lra & \V_v(\T P) \\
u & \lmt & \left.\frac{\d}{\d t}\right|_{t=0}v+t u\,.
\eenn
Let $\mu$ be in $\bigwedge^2\T_x P$.
The connection also defines a splitting of the tangent space of $\bigwedge^2\T P$ at $\mu$ into
the direct sum of a horizontal space isomorphic to $\T_x P$ and a vertical space isomorphic
to $\bigwedge^2\T_x P$. Because $p\circ\eta$ is the identity of $P$, the horizontal component
of $T\eta:\T P\rightarrow \T\bigwedge^2\T P$ is the identity. Denote its vertical component by
$$\D\eta:\T P\lra\bigwedge^2\T P.$$
If $f$ is a smooth function on $\T P$, its differential at $v$ in $T_x P$ has a horizontal
and a vertical component. Denote by $\del_2f:T_x P\lra\C$ its horizontal component and
by $\del_1f:T_x P\lra\C$
its vertical one. In the same fashion, for $\pi:\T P\rightarrow P$, denote by $\del_1\pi$ and
$\del_2p$ respectively the vertical and horizontal components of
$\pi_*:\T(\T P)\rightarrow \T P$.

\begin{lemm}\label{lemm:bracketonTP}
Let $f$ and $g$ be smooth functions on $\T P$. For $x$ in $P$ and $v$ in $T_x P$
$$\{f,g\}_{\T P}(v)=\lb\D\eta(v),\del_1f(v)\wedge\del_1g(v)\rb+
\lb\eta(x),\del_1f(v)\wedge\del_2g(v)-\del_1g(v)\wedge\del_2f(v)\rb.$$
\end{lemm}
\begin{proof}
Use the right hand side of the above equation to define a bracket
\benn
\{\,,\}:\;\C^\infty(P)\times\C^\infty(P) & \lra & \C^\infty(P).
\eenn
This bracket satisfies the Leibniz identity because the operators $\del_1$ and $\del_2$
do. It is also clearly anti-symmetric. To prove that it is equal to $\{\,,\}_{\T P}$,
it suffices to prove that it satisfies
Equations~(\ref{eqn:bracketonTV1}), (\ref{eqn:bracketonTV2}) and (\ref{eqn:bracketonTV3}).

Equation~(\ref{eqn:bracketonTV2}) is satisfied because $\del_1$ vanishes on pull-backs
to $\T P$ of functions on $P$.

If $g$ is a function on $P$, then $\del_2 p^*g=\d g$. If $f$ is equal to $\widetilde{\alpha}$
for some $1$-form $\alpha$ on $P$, then $\partial_1\widetilde{\alpha}(v)=\alpha_x$ for all
$v$ in $T_x P$.
Hence the bracket $\{\,,\}$ satisfies
Equation~(\ref{eqn:bracketonTV3}).

The connection on $\T P\rightarrow P$ also defines a connection on its dual bundle and on
all bundles one can construct from $\T P$ and $\T^*P$ through direct sums, tensor products... 

Let $\alpha$ and $\beta$ be closed $1$-forms on $P$. Let $v$ be in $\T_x P$ and let $\sigma$
be a path in $P$ such that $\sigma(0)=x$ and $\dot{\sigma}(0)=v$; for example, take
$\sigma(t)=\exp(tv)$, where $\exp$ is the exponential map of the connection.
Equation~(\ref{eqn:bracketonTV1}) gives
\be
\{\widetilde{\alpha},\widetilde{\beta}\}_{\T P}(v)=
\left.\frac{\d}{\d t}\right|_{t=0}\eta_{\sigma(t)}(\alpha\wedge\beta).\label{eqn:bracketonTP}
\ee
Let $E\rightarrow P$ be the bundle $\T^*P\oplus\T^*P\oplus\bigwedge^2\T P\rightarrow P$. There is a
natural map
\benn
m:\,E & \lra & \mathbb{R} \\
(\alpha,\beta,\mu) & \lmt & \mu(\alpha\wedge\beta)
\eenn
I will compute the differential of this map. Let $w$ be in $\T_{(\alpha_x,\beta_x,\mu_x)}E$
with horizontal component $\H(w)=v$.
Its vertical component $\V(w)=(\epsilon_1,\epsilon_2,\theta)$ is in
$\T^*_x P\oplus\T^*_x P\oplus\wedge^2\T_x P$. Let
$\phi_{\sigma(t)}:E_x\rightarrow E_{\sigma(t)}$ be the parallel
transport along the path $\sigma$. Define a path in $E$ by
$$\gamma(t)=\phi_{\sigma(t)}(\alpha_x+t\epsilon_1,\beta_x+t\epsilon_2,\mu_x+t\theta).$$
The path $\gamma$ satisfies
$$\gamma(0)=(\alpha_x,\beta_x,\mu_x)\mbox{ and }\dot{\gamma}(0)=w.$$
Notice that, because the connection on $E$ is defined using a single connection on
$\T P$, we have
$$m\circ\phi_{\sigma(t)}=m.$$
This means
\benn
m_*(w) & = & \left.\frac{\d}{\d t}\right|_{t=0}m\circ\phi_{\sigma(t)}
(\alpha_x+t\epsilon_1,\beta_x+t\epsilon_2,\mu_x+t\theta) \\
 & = & m(\epsilon_1, \beta_x,\mu_x)+m(\alpha_x,\epsilon_2,\mu_x)+m(\alpha_x,\beta_x,\theta).
\eenn
This last computation together with Equation~(\ref{eqn:bracketonTP}) gives
\be
& & \{\widetilde{\alpha},\widetilde{\beta}\}_{T P}(v)=m(\partial_{1,x}\alpha(v),\beta_x,\eta_x)+
m(\alpha_x,\partial_{1,x}\beta(v),\eta_x)+m(\alpha_x,\beta_x,\D\eta(v)).\hspace*{-2cm}\label{equa:intermediary-bracket-on-TP}
\ee
Firstly, in this equality, one can replace $\alpha_x$ and $\beta_x$ by respectively
$\partial_{1,v}\widetilde{\alpha}$ and $\partial_{1,v}\widetilde{\beta}$.\\
Secondly, consider a vector field $X$ on $P$. Let
\benn
\iota:P & \lra & T^*P\oplus T P\\
x & \lmt & (\alpha_x,X_x),
\eenn
and
\benn
k:T^*P\oplus T P & \lra & \C\\
(\delta,Z) & \lmt & \delta(Z).
\eenn
I choose $X$ such that $\partial_1 X=0$. Then, differentiation of the equality
\benn
k\circ\iota=\widetilde{\alpha}\circ X
\eenn
leads to
\be
\partial_{2,X_x}\widetilde{\alpha}(v)=\partial_{1,x}\alpha(v)(X_x).\label{equa:small-equation}
\ee
Since $\alpha$ and $\beta$ are closed $1$-forms, and since $\nabla$ is torsion free
\be
\partial_{1,x}\alpha(v)(X_x)=\partial_{1,x}\alpha(X_x)(v).\label{equa:2small-equation}
\ee
Equations~(\ref{equa:intermediary-bracket-on-TP}),~(\ref{equa:small-equation}) and
(\ref{equa:2small-equation}) put together prove that the bracket $\{\,,\}$ satisfies
Equation~(\ref{eqn:bracketonTV1}).
\end{proof}

I deduce the following corollaries.
\begin{coro}\label{coro:poissonmapforriemanniancompatiblestructure}
Let $\pi$ be a map $\T P\rightarrow P$. Let $f$ and $g$ be functions on $P$. Their pull-backs by $\pi$ satisfy
\benn
\hspace*{-1cm}\{\pi^*f,\pi^*g\}_{\T P}(v) & = & \lb\,(\del_1\pi(v)\otimes\del_1\pi(v))(\D_x\eta(v))+
2(\del_1\pi(v)\otimes\del_2\pi(v))(\eta(x)),\d_{\pi(v)}f\wedge d_{\pi(v)}g\,\rb.
\eenn
In particular, $\pi$ is Poisson if and only if
\benn
\frac{1}{2}(\del_1\pi(v)\odot\del_1\pi(v))(\D_x\eta(v))+(\del_1\pi(v)\odot\del_2\pi(v))(\eta(x)) & = & \eta(\pi(v)),
\eenn
for all $v$ in $\T_x P$, where $\odot$ means the symmetric product.
\end{coro}
\begin{coro}
Assume the Poisson bivector field is parallel relative to the connection $\nabla$.
Then, a map $\pi:\T P\rightarrow P$ is Poisson if and only if
\benn
(\del_1\pi(v)\odot\del_2\pi(v))(\eta(x)) & = & \eta(\pi(v)),
\eenn
for all $v$ in $\T_x P$.
\end{coro}
\demo
Indeed, since the Poisson bivector field is parallel
$$\D\eta(v)=0.$$
\findemo
Assume the connection on $P$ is the Levi-Civita connection of a metric on $P$ and assume $P$ is complete.
Its exponential map is denoted
$$\exp:\T P\rightarrow P.$$
When restricted to a fibre $\T_x P$, I will denote it $\exp_x$. I quote here the following
Lemma for future reference.
\begin{lemm}
Let $v$ be a tangent vector to $P$ at a point $x$. Let $w$ be a tangent vector to $\T P$ at $v$. Its horizontal and
vertical components are respectively $\H w$ and $\V w$. Consider the geodesic $\sigma(t)=\exp(t\H w)$ and the
$1$-parameter family of geodesics
$$\gamma_s(t)=\exp_{\sigma(t)}(s\phi_{\sigma(t)}(v+t\V w)).$$
The differential of $\exp$ at $v$ is given by
$$\T_v\exp(w)=\left.\frac{\d}{\d t}\right|_{t=0}\gamma_1(t).$$
It follows that $\T_v\exp(w)$ is also the value at $t=1$ of the Jacobi field $J$ along the geodesic $\sigma$ with
initial value $J(0)=\H w$ and $J^\prime(0)=\V w$. In particular, $\exp:\T P\rightarrow P$ is a submersion.
\end{lemm}
\demo
This is a simple exercise in Riemannian geometry.
\findemo


\section{The torus with constant Poisson structure}\label{sect:the-torus}
Let $P$ be the $n$-dimensional torus $\R^n/2\pi\Z^n$ with its metric inherited
from the euclidean metric $\lb\,,\rb$ on $\R^n$. Consider a constant
Poisson structure on $P$ given by a skew-symmetric $n\times n$ matrix $\eta$.
Identify $\T P$ with $\R^n\times P$ in
the obvious way. From
Corollary~\ref{coro:poissonmapforriemanniancompatiblestructure}, I deduce that the map
\benn
\pi:\R^n\times P & \lra & P\\
(u,p) & \lmt & \exp_p(\frac{1}{2}u)=p+\frac{1}{2}u
\eenn
is Poisson. It is also a surjective submersion. Hence , I can hope that a quantisation
of $\T P$ will lead to a quantisation of $P$.

The dual $A=T^*P$ of $\T P$, identified with $\T P$ using the euclidean metric, is a Lie algebroid.
It can be integrated to a source simply-connected Lie groupoid. The space of
morphisms of this groupoid is $\R^n\times P$. Notice that $\R^n$ is the direct orthogonal
sum of $\ker(\eta)$ and $\im(\eta)$. Let $\pr_1$ be the orthogonal projection on $\ker(\eta)$
and $\pr_2$ the orthogonal projection on $\im(\eta)$. The source map of the groupoid is
\benn
s:\R^n\times P & \lra & P\\
(u,p) & \lmt & p,
\eenn
whereas the target map is
\benn
t:\R^n\times P & \lra & P\\
(u,p) & \lmt & p+\pr_2(u).
\eenn
Given two elements $(u,p)$ and $(v,q)$ in the groupoid, their multiplication $(u,p)\cdot(v,q)$ is well-defined
if the target of $(u,p)$ is equal to the source of $(v,q)$, that is if $q=p+\pr_2(u)$; in this case
$$(u,p)\cdot(v,q)=(u+v,p).$$

Assume that $\eta$ is invertible, that is $P$ is symplectic. In this case, $\eta$ defines
an isomorphism of Lie algebroids between $T^*P$ and $T P$
\benn
T^*P & \lra & T P\\
(\xi,p) & \lmt & (\eta(\xi),p).
\eenn
Choose the natural connection on the trivial
vector bundle $\T P\simeq\R^n\times P\lra P$.
The exponential map $\exp$ for $T P$ is
\benn
\T P & \lra & \R^n\times P \\
(X,p) & \lmt & (X,p).
\eenn
Whereas the $\exp$ map for $T^*P$, that is $\alpha$, is
\benn
\T^*P & \lra & P \\
(\xi,p) & \lmt & (\eta(\xi),p).
\eenn
The tangent groupoid is given by $\widetilde{G}=\R\times\R^n\times P$ with
\benn
s(\hbar,u,q)=(\hbar,q)\mbox{, and }t(\hbar,u,q)=(\hbar,q+\hbar u).
\eenn
The product is given by
\benn
(\hbar,v,q+\hbar u)\cdot(\hbar,u,q)=(\hbar,v+u,q).
\eenn
With this representation of $\widetilde{G}$, the exponential of the tangent groupoid is
\benn
\R\times T P & \lra & \R\times\R^n\times P\\
(\hbar,X,q) & \lmt & (\hbar,X,q),
\eenn
therefore $\widetilde{\alpha}$ is
\benn
\R\times T^*P & \lra & \R\times\R^n\times P\\
(\hbar,\xi,q) & \lmt & (\hbar,\eta(\xi),q).
\eenn
Let $f$ be a function on $A^*=T P$. Assume it is acceptable for quantisation.
Let $H$ be a compactly supported smooth function on $\widetilde{G}$. For
$(\hbar,q)$ in $\R\times P$,
\be
\mathcal{Q}(f)_{\hbar,q}(H)=\int_{\R^n}\d(X)\, f(X,q)\int_{\R^n}\d(\xi)
e^{-i\lb\xi,X\rb}H(\hbar,\eta(\xi),q).\label{equa:Q(f)-on-P}
\ee
Let $r$ be a vector in $\Z^n$ and define a function
\benn
g_r:P & \lra & \C\\
q & \lmt & e^{i\lb r,q\rb}.
\eenn
The number $e^{i\lb r,q\rb}$ is well defined because $\lb r,q\rb$ is well defined
modulo $2\pi$. Set $f_r=\pi^*g_r$, that is
\benn
f_r:T P & \lra & \C\\
(X,q) & \lmt & e^{i\lb r,q\rb}e^{\half i\lb r,X\rb}.
\eenn
\begin{prop}\label{prop:Q(fr)-and-Q(fs)}
Let $r$ and $r^\prime$ be vectors in $\Z^n$, then
\benn
\mathcal{Q}(f_r)\star\mathcal{Q}(f_{r^\prime})=
\mathcal{Q}(e^{i\frac{\hbar}{2}\lb r,\eta(s)\rb}f_{r+r^\prime}).
\eenn
\end{prop}
\begin{proof}
With $f_r$ instead of $f$, (\ref{equa:Q(f)-on-P}) becomes
\benn
\mathcal{Q}(f_r)_{\hbar,q}(H) & = & \int_{\R^n}\d(X)\, e^{i\lb r,q\rb}e^{\half i\lb r,X\rb}
\int_{\R^n}\d(\xi)e^{-i\lb\xi,X\rb}H(\hbar,\eta(\xi),q).\label{equa:Q(f)-on-P}\\
 & = & e^{i\lb r,q\rb}H(\hbar,\eta(\half r),q).
\eenn
Hence, the product $\mathcal{Q}(f_r)\star\mathcal{Q}(f_{r^\prime})$ is
\benn
\hspace*{-1cm}\mathcal{Q}(f_r)\star\mathcal{Q}(f_{r^\prime})_{\hbar,q}(H) & = &
\int_{\R^n}\d(X)\;f_{r^\prime}(X,q)\int_{\R^n}\d(\xi)\;e^{i\lb\xi,X\rb}\int_{\R^n}\d(Y)\;
f_r(Y,q+\hbar\eta(\xi))\\
 & & \int_{\R^n}\d(\zeta)\;e^{-i\lb\zeta,Y\rb}H(\hbar,\eta(\zeta+\xi),q)\\
 & = & \int_{\R^n}\d(X)\;e^{i\lb r^\prime ,q\rb}e^{\half i\lb r^\prime ,X\rb}\int_{\R^n}\d(\xi)\;e^{i\lb\xi,X\rb}\int_{\R^n}\d(Y)\;
  e^{i\lb r,q+\hbar\eta(\xi)\rb}e^{\half i\lb r,Y\rb}\\
 & & \int_{\R^n}\d(\zeta)\;e^{-i\lb\zeta,Y\rb}H(\hbar,\eta(\zeta+\xi),q)\\
 & = & e^{i\lb r^\prime +r,q\rb}\int_{\R^n}\d(Y)\;e^{i\lb r,\hbar\eta(\half r^\prime )\rb}e^{i\lb r,Y\rb}
   \int_{\R^n}\d(\zeta)\; e^{-i\lb\zeta,Y\rb}H(\hbar,\eta(\zeta+\half r^\prime ),q)\\
 & = & e^{i\lb r^\prime +r,q\rb}e^{i\frac{\hbar}{2}\lb r,\eta(r^\prime )\rb}H(\hbar,\half\eta(r+r^\prime ),q)\\
 & = & e^{i\frac{\hbar}{2}\lb r,\eta(r^\prime )\rb}\mathcal{Q}(f_{r+r^\prime })_{\hbar,q}(H).
\eenn
\end{proof}
Let $\mathcal{P}$ be the algebra of functions on $P$ generated by $\{g_r,\,r\in\Z^n\}$. It
is a dense sub-algebra of the $C^*$-algebra of continuous functions on $P$.
Proposition~\ref{prop:Q(fr)-and-Q(fs)} shows that the product on this algebra can be deformed in
\benn
g_r\star_\hbar g_{r^\prime}=e^{i\frac{\hbar}{2}\lb r,\eta(s)\rb}g_{r+r^\prime},
\eenn
for each $\hbar$. With this new product, $\mathcal{P}$ becomes a $*$-algebra which can be
completed into a $C^*$-algebra $\mathcal{P}_\hbar$. The natural family of injections
$\mathcal{P}\lra\mathcal{P}_\hbar$ gives the usual quantisation of the torus with constant
Poisson structure as defined in Tang and Weinstein~\cite{TW}. It is
a strict deformation quantisation in the
sense of Rieffel.


\section*{Appendix: The $2$-sphere in $\R^3$}
Let $P=S^2$ be the $2$-sphere $\{(x,y,z)\in\mathbb{R}^3/x^2+y^2+z^2=1\}$.
In this Section, I show how to construct a Poisson map between $T S^2$ and $S^2$.

Consider the metric on $S^2$ given by the restriction of the euclidean metric
$\d x^2+\d y^2+\d z^2$ on $\mathbb{R}^3$. For $p=(x,y,z)$
in $\mathbb{R}^3$, define the endomorphism
$$J_p(u)=u\wedge p\mbox{, for all $u\in\mathbb{R}^3$.}$$
The restriction of $J$ to each tangent space of $S^2$ defines a complex structure
on the sphere. Also
$$\omega=g(J\cdot,\cdot)$$
is a symplectic form on $S^2$. It is the restriction to $S^2$ of the $2$-form
$z\d y\wedge\d z-y\d x\wedge\d z+z\d x\wedge\d y$ defined on $\mathbb{R}^3$. The
geodesics are great circles on the sphere so that the exponential is given by
$$\exp_p(u)=p\cos(\norm{u})+\frac{\sin(\norm{u})}{\norm{u}}u.$$
The differential $\del_1\exp$ is a map from $T_p S^2$ to $T_{\exp_(u)}S^2$ given by
\benn
\del_1\exp_{(u,p)}(h) & = & \left.\frac{\d}{\d t}\right|_{t=0}\exp_p(u+th) \\
 & = & \frac{\sin(\norm{u})}{\norm{u}}(-\lb u,h\rb p+h)+\frac{\norm{u}\cos(\norm{u})-\sin(\norm{u})}{\norm{u}^2}
\lb\frac{u}{\norm{u}},h\rb u
\eenn

In addition,
$$\del_2\exp_{(u,p)}(\epsilon) = \left.\frac{\d}{\d t}\right|_{t=0}\exp_{\sigma(t)}\phi_{\sigma(t)}(u),$$
where $\sigma(t)=\exp_p(t\epsilon)$ and $\phi_{\sigma(t)}$ is the parallel
transport along $\sigma(t)$. Without loss
of generality, I can assume that $\epsilon$ is a unit vector. In this situation,
$$\sigma(t)=\cos(t)p+\sin(t)\epsilon.$$
Also, $\{p,\epsilon,p\wedge\epsilon\}$ forms an orthonormal basis of $\mathbb{R}$
and $\phi_{\sigma(t)}$ is a morphism in $\mbox{SO}(3)$. It is given by
\benn
\phi_{\sigma(t)}(p) & = & \sigma(t) \\
\phi_{\sigma(t)}(\epsilon) & = & \frac{\d}{\d t}\sigma(t)\\
 & = & -\sin(t)p+\cos(t)\epsilon\\
\phi_{\sigma(t)}(p\wedge\epsilon) & = & p\wedge\epsilon.
\eenn
It follows that
$$\phi_{\sigma(t)}(u) = \lb u,\epsilon\rb(-\sin(t)p+\cos(t)\epsilon)+\lb u,p\wedge\epsilon\rb p\wedge\epsilon.$$
The parallel transport preserves the norm, hence $u$ and $\phi_{\sigma(t)}(u)$
have the same norm and
\benn
\exp_{\sigma(t)}\phi_{\sigma(t)}(u) & = & \cos(\norm{u})\sigma(t)
  +\frac{\sin(\norm{u})}{\norm{u}}\phi_{\sigma(t)}(u)\\
 & = & \cos(\norm{u})(\cos(t)p+\sin(t)\epsilon)+\\
 & & +\frac{\sin(\norm{u})}{\norm{u}}(\lb u,\epsilon\rb(-\sin(t)p+\cos(t)\epsilon)
  +\lb u,p\wedge\epsilon\rb p\wedge\epsilon).
\eenn
I can now compute
\benn
\del_2\exp_{(u,p)}(\epsilon) & = & \left.\frac{\d}{\d t}\right|_{t=0}\exp_{\sigma(t)}\phi_{\sigma(t)}(u)\\
 & = & \cos(\norm{u})\epsilon-\frac{\sin(\norm{u})}{\norm{u}}\lb u,\epsilon\rb p.
\eenn
This last formula is of course still valid when $\epsilon$ is not a unit vector.

Consider a map
\benn
\pi:\T S^2 & \lra & S^2\\
u & \lmt & \ba{l}\exp(\lambda\norm{u}) \\
=\cos(\lambda\norm{u})p+f(\lambda\norm{u})\lambda u,\ea
\eenn
where $\lambda$ is a function of $\norm{u}^2$ defined for $\norm{u}^2$ in
some neighbourhood of $0$ in $\R$. For $\pi$ to be a Poisson map, I need
$\lambda(0)=\half$.

I will now compute the differential of such a map. Firstly,
\def\dif{{\left.\frac{\d}{\d t}\right|_{t=0}}}
\def\fl{{f(\lambda\norm{u})}}
\def\fpl{{f^\prime(\lambda\norm{u})}}
\def\uh{{\lb u,h\rb}}
\def\lu{{\lambda\norm{u}}}
\benn
\partial_{1,(u,p)}\pi(h) & = &
\dif\exp_p(\lambda(\norm{u+t h}^2)(u+t h))\\
& = & \dif\cos(\lambda(\norm{u+t h}^2)\norm{u+t h})p+\\
& & +f(\lambda(\norm{u+t h}^2)\norm{u+t h})\lambda(\norm{u+t h}^2)(u+t h)\\
& = & -\sin(\lambda\norm{u})(2\lambda^\prime\lb u,h\rb\norm{u}+
 \lambda\lb h,\frac{u}{\norm{u}}\rb)p+\fl\lambda h+\\
& & +2\fl\uh\lambda^\prime u+(2\uh\lambda^\prime\norm{u}
 +\lambda\lb h,\frac{u}{\norm{u}}\rb)\fpl\lambda u\\
& = & -\lambda\fl(2\lambda^\prime\uh\norm{u}^2+\lambda\uh)\lambda\fl p+\fl\lambda h+\\
& & +(2\fl\lambda^\prime+2\lambda^\prime\lambda\fpl+\frac{\lambda^2}{\norm{u}}\fpl)\uh u\\
& = & -\lambda\fl(2\lambda^\prime\uh\norm{u}^2+\lambda\uh)\lambda\fl p+\fl\lambda h+\\
& & +(2\lambda^\prime\cos(\lambda\norm{u})+\frac{\lambda^2}{\norm{u}}\fpl)\uh u,
\eenn
where for the last equality I have used the relation
$$tf^\prime(t)=\cos(t)-f(t).$$
Secondly, assuming without loss of generality that $\norm{\epsilon}=1$,
\benn
\partial_{2,(u,p)}\pi(\epsilon) & = &
\dif\pi(\phi_{\sigma(t)}(u))\\
& = & \dif\exp_{\sigma(t)}(\lambda(\norm{\phi_{\sigma(t)}(u)}^2)\phi_{\sigma(t)}(u))\\
& = & \dif\exp_{\sigma(t)}(\lambda(\norm{(u)}^2)\phi_{\sigma(t)}(u))\\
& = & \dif\cos(\lambda\norm{\phi_{\sigma(t)}(u)})\sigma(t)
  +f(\lambda\norm{\phi_{\sigma(t)}(u)})\lambda\phi_{\sigma(t)}(u)\\
& = & \dif\cos(\lu)\sigma(t)+\fl\lambda\phi_{\sigma(t)}(\lb u,\epsilon\rb\epsilon
  +\lb u,p\times\epsilon\rb p\times\epsilon)\\
& = & \cos(\lu)\sigma^\prime(0)+\lambda\fl\dif(\lb u,\epsilon\rb\sigma^\prime(t)
  +\lb u,p\times\epsilon\rb p\times\epsilon)\\
& = & \cos(\lu)\epsilon+\lb u,\epsilon\rb\lambda\fl\sigma^{\prime\prime}(0)\\
& = & \cos(\lu)\epsilon-\lb u,\epsilon\rb\lambda\fl p.
\eenn
I wish to compute $\norm{u}^2\partial_1\pi\odot\partial_2\pi(\eta(p))$. I know
$$\eta(p)=\frac{1}{\norm{u}^2}u\wedge(p\times u)\mbox{, whenever $u\neq 0$.}$$
So, I need to compute
\benn
\partial_{1,(u,p)}\pi(u) & = & -\lambda\fl(2\lambda^\prime\norm{u}^2+\lambda)\norm{u}^2p+\lambda\fl u+\\
 & & +(2\lambda^\prime\cos(\lu)+\frac{\lambda}{\norm{u}^2}\fpl)\norm{u}^2u\\
& = & -\lambda\norm{u}^2\fl(2\lambda^\prime\norm{u}^2+\lambda)p+(\lambda\fl+\\
 & & +2\lambda^\prime\cos(\lu)\norm{u}^2+\lambda^2\norm{u}\fpl)u\\
& = & (2\lambda^\prime\norm{u}^2+\lambda)(\cos(\lu)u-\norm{u}^2\lambda\fl p),
\eenn
and
\benn
\partial_{1,(u,p)}\pi(p\times u) & = & \lambda\fl p\times u,
\eenn
and
\benn
\partial_{2,(u,p)}\pi(u) & = & \cos(\lu)u-\norm{u}^2\lambda\fl p,
\eenn
and finally
\benn
\partial_{2,(u,p)}\pi(p\times u) & = & \cos(\lu)p\times u.
\eenn
Notice that
$$\partial_{1,(u,p)}\pi(u) = (2\lambda^\prime\norm{u}^2+\lambda)\partial_{2,(u,p)}\pi(u).$$
I can now compute
\benn
\norm{u}^2\partial_1\pi\odot\partial_2\pi(\eta(p)) & = & 
  \partial_1\pi(u)\wedge\partial_2\pi(p\times u)+\partial_2\pi(u)\wedge\partial_2\pi(p\times u)\\
& = & \partial_2(u)\wedge((2\lambda^\prime\norm{u}^2+\lambda)\partial_{2,(u,p)}\pi(u)+\partial_1\pi(p\times u))\\
& = & (\cos(\lu)u-\norm{u}^2\lambda\fl p)\wedge((2\lambda^\prime\norm{u}^2+\lambda)\cos(\lu)+\\
 & & +\lambda\fl)p\times u\\
& = & (\cos(\lu)u-\norm{u}^2\sin(\lu)p)\wedge((2\lambda^\prime\norm{u}^2+\lambda)\cos(\lu)+\\
 & & +\lambda\fl)p\times u.
\eenn
On the other hand
\benn
\norm{u}^2\eta(\pi(u,p)) & = & \norm{u}^2\phi_{\sigma(\lu)}(\eta(p))
  \;\;\;\;\;\;\;\mbox{, with }\sigma(t) = \exp(t\frac{u}{\norm{u}})= \cos(t)p+\sin(t)\frac{u}{\norm{u}}\\
& = & \phi_{\sigma(\lu)}(u)\wedge\phi_{\sigma(\lu)}(p\times u)\\
& = & \sigma^\prime(\lu)\wedge(p\times u)\\
& = & (\cos(\lu)u-\norm{u}\sin(\lu)p)\wedge(p\times u).
\eenn
It follows that $\pi$ is a Poisson map if and only if $\lambda$ satisfies the following
differential equation
$$(2\lambda^\prime(t^2)t^2+\lambda(t^2))\cos(\lambda(t^2)t)+\lambda(t^2)f(\lambda(t^2)t)=1.$$
Put $\mu(t)=\lambda(t^2)$ so that $\mu^\prime(t)=2t\lambda^\prime(t^2)$. The function
$\mu$ satisfies the differential equation
$$(\mu^\prime(t)t+\mu(t))\cos(\mu(t)t)+\frac{\sin(\mu(t)t)}{t}=1.$$
Put $\alpha(t)=\sin(\mu(t)t)$ so that $\alpha^\prime(t)=(\mu^\prime(t)t+\mu(t))\cos(\mu(t)t)$.
The function $\alpha$ satisfies the differential equation
$$t\alpha^\prime(t)+\alpha(t)=t.$$
A general solution of this equation is
$$\alpha(t)=\frac{a}{t}+\frac{t}{2}\mbox{, with }a\in\R.$$
Hence
$$\mu(t)=\frac{1}{t}\arcsin(\frac{a}{t}+\frac{t}{2}).$$
Since I want $\mu(0)=\lambda(0)=\half$, I need $a=0$ and
$$\mu(t)=\frac{1}{t}\arcsin(\frac{t}{2}).$$
I deduce
\begin{prop}
The map
\benn
\pi:\T S^2 & \lra & S^2\\
(u,p) & \lmt & \exp_p(\frac{1}{\norm{u}}\arcsin(\frac{\norm{u}}{2})u)
\eenn
is Poisson
\end{prop}
The way it is written in the previous Proposition, the map $\pi$ is only defined on
a neighbourhood of $S^2$ in $T S^2$. Nevertheless, it can easily be extended to a
continuous function on the whole of $T S^2$, so that it is Poisson wherever it is
smooth. The fact that this map is not smooth at all points means that technics used
in Section~\ref{sect:the-torus} will not carry over here. I nevertheless believe that
a modification of these technics will make things work.


{\providecommand{\bysame}{\leavevmode ---\ }}


\begin{thebibliography}{1}

%
\bibitem{Co}{
\scshape A.~Connes},
-- \emph{Noncommutative geometry}, Academic Press, 1994.
%
\bibitem{HS}{
\scshape M.~Hilsum\normalfont{ and }G.~Skandalis},
-- \emph{Morphismes $K$-orient\'es d'espaces de feuilles et fonctorialit\'e
en th\'eorie de Kasparov (d'apr\`es une conjecture d'A.~Connes). (French)
[$K$-oriented morphisms of spaces of leaves and functoriality in Kasparov theory
(after a conjecture of A. Connes)]} Ann. Sci. \'Ecole Norm. Sup. (4) 20 (1987), no. 3, 325--390.
%
\bibitem{La}{
\scshape N.~P.~Landsman},
-- \emph{Mathematical topics between classical and quantum mechanics}, Springer.
%
\bibitem{LR}{
\scshape N.~P.~Landsman, B.~Ramazan},
-- \emph{Quantization of Poisson algebras associated to Lie algebroids},
Ramsay, Arlan (ed.) et al., Groupoids in analysis, geometry, and physics.
AMS-IMS-SIAM joint summer research conference, University of Colorado,
Boulder, CO, USA, June 20-24, 1999. Providence, RI: American Mathematical
Society (AMS). Contemp. Math. 282, 159-192 (2001).
%
\bibitem{NWX}{
\scshape V.~Nistor, A.~Weinstein\normalfont{ and }P.~Xu},
-- \emph{Pseudodifferential operators on differential groupoids},
Pacific J. Math. 189 (1999), no. 1, 117--152.
%
\bibitem{Ra}{
\scshape B.~Ramazan},
-- \emph{Quantification par d\'eformation des vari\'et\'es de Lie-Poisson}, PhD Thesis,
Universit\'e d'Orl\'eans (1998).
%
\bibitem{Rie}{
\scshape M.~.A~Rieffel},
-- \emph{Deformation quantization of Heisenberg manifolds},
Commun. Math. Phys. \bf{122} \rm (1989), 531--562.
%
\bibitem{TW}{
\scshape X.~Tang\normalfont{ and }A.~Weinstein},
-- \emph{Quantization and Morita Equivalence for Constant Dirac Structures on Tori},
arXiv:math.QA/0305413 (2003).
%
\bibitem{Wei}{
\scshape A.~Weinstein},
-- \emph{Symplectic groupoids, geometric quantization, and irrational rotation algebras},
Symplectic geometry, groupoids, and integrable systems, S\'eminaire sud-Rhodanien de g\'eom\'etrie
\`a Berkeley (1989), P.~Dazor and A.~Weinstein, eds., Springer-MSRI Series (1991), 281--290.
%
\bibitem{Wei2}{
\scshape A.~Weinstein},
-- \emph{Blowing up realizations of Heisenberg-Poisson manifolds}, Bull. Sci. Math. 113 (1989),
no. 4, 381--406.
%
\bibitem{Wei3}{
\scshape A.~Weinstein},
-- \emph{Tangential deformation quantization and polarized symplectic groupoids},
Deformation theory and symplectic geometry (Ascona, 1996), 301--314,
Math. Phys. Stud., 20, Kluwer Acad. Publ., Dordrecht, 1997.
%
\end{thebibliography}
\end{document}